\documentclass[11pt,leqno,english]{article}

\usepackage{646}

\begin{document}



\title{\protect
	On the number of ${L}_{\infty\omega_1}$-equivalent
	non-isomorphic models
}

\author{
	Saharon Shelah
	\thanks{Thanks to GIF for its support of this research and
	also to University of Helsinki for funding a visit of the
	first author to Helsinki in August 1996. Pub. No. 646.} \\
	\and
	Pauli V\"{a}is\"{a}nen
	\thanks{This paper is the second author's Licentiate's
	thesis. The second author did his share of the paper under the
	supervision of Tapani Hyttinen.}
}

\date{\today}

\maketitle

\begin{abstract}%
 We prove that if \ZF is consistent then \ZFC + \GCH is consistent
with the following statement: There is for every $k < \omega$ a model
of cardinality $\aleph_1$ which is \Lan-equivalent to exactly $k$
non-isomorphic models of cardinality $\aleph_1$. In order to get this
result we introduce ladder systems and colourings different from the
``standard'' counterparts, and prove the following purely
combinatorial result: For each prime number $p$ and positive integer
$m$ it is consistent with \ZFC + \GCH that there is a ``good'' ladder
system having exactly $p^m$ pairwise nonequivalent colourings.%
 \footnote{%
}%
\end{abstract}


\begin{SECTION} {-} {Introduction} {Introduction}




If \M is a model, \Card \M denotes the cardinality of the universe of
\M. Suppose \M and \N are two models of the same vocabulary and
$\kappa$ is a cardinal. We write $\M \LEquiv [\kappa] \N$ if \M and \N
satisfy the same sentences of the infinitary language \Lan
[\kappa]. For a definition of \Lan [\kappa], the reader is referred to
\cite {Dickmann}. For any model \M of cardinality $\kappa$, define
 \[
    \No \M = \Card [\Big] {
	\Set [\big] {\Quotient \N \Isomorphic} {
		\Card \N = \kappa \And
		\N \LEquiv [\kappa] \M 
	} 
    }, 
 \]
where \Quotient \N \cong is the equivalence class of \N under the
isomorphism relation. We study the possible values of \No \M for
models \M of cardinality $\aleph_1$. In particular, we prove the
following theorem:

\begin{THEOREM}{Models}
 Assuming \ZF is consistent, it is consistent with \ZFC + \GCH that
there is for every $k < \omega$ a model \M \Note {of a vocabulary of
cardinality $\leq \aleph_1$} such that $\Card \M = \aleph_1$ and $\No
\M = k$.
 \end{THEOREM}

When \M is countable, $\No \M = 1$ by \cite {Scott}. This result
extends to structures of cardinality $\kappa$ when $\kappa$ is a
singular cardinal of countable cofinality \cite {Chang}. So the study
of possible values of \No \M is divided into the following cases
according to the cardinality of \M:
\begin{enumerate}

\ITEM{Weakly}%
\Card \M is weakly compact;

\ITEM{Singular}%
\Card \M is singular of uncountable cofinality;

\ITEM{Regular}%
\Card \M is uncountable, regular, and non-weakly compact.

\end{enumerate}

In \cite {Sh133} Shelah was able to show that when $\kappa$ is a
weakly compact cardinal there is for every non-zero cardinal $\mu \leq
\kappa$, a model $\M$ such that $\Card \M = \kappa$ and $\No \M =
\mu$. In a paper which is in preparation by the authors, the problem
of the possible value of \No \M between $\kappa$ and $2^\kappa$ for a
model \M of weakly compact cardinality is completely solved.

Shelah has considered the singular case in two of his papers \cite
{Sh189, Sh228}. Let $\kappa$ be a singular cardinal of uncountable
cofinality. In the former paper it is shown that if one allows
relation symbols of arbitrary large arity $< \kappa$ and $\mu$ is a
non-zero cardinal with $\mu^{\Cf \kappa} < \kappa$, then there exists
a model \M of singular cardinality $\kappa$ with $\No \M = \mu$. In
the latter paper Shelah gives a general way to build models \M with
relations of finite arity only and for which the value of \No \M is
quite arbitrary: for every non-zero cardinal $\mu \in \kappa \Union
\Braces {\kappa^{\Cf \kappa}}$, there exists a model \M of cardinality
$\kappa$ such that $\No \M = \mu$ and its vocabulary consists of one
binary relation symbol, provided that $\theta^{\Cf \kappa} < \kappa$
for all $\theta < \kappa$. The paper \cite {Sh228} together with a
recent paper \cite {ShVa644} offer a complete answer to the singular
case provided that the singular cardinal hypothesis holds. For example
it follows that $\No \M = \kappa$ is possible, even in $L$.


If $V = L$ and $\kappa \geq \aleph_1$ is a regular cardinal which is
not weakly compact, \No \M has either the value 1 or $2^\kappa$ for
all models \M having cardinality $\kappa$. For $\kappa = \aleph_1$
this result was first proved in \cite {Palyutin}. Later Shelah
extended the result to all other regular non-weakly compact cardinals
in \cite {Sh129}.

It seems that there are no published independence results about the
case that \Card \M is a regular but not weakly compact cardinal. But
it is known that the independence result given in \cite {Sh125}
implies the consistency of ``there is a model \M of cardinality
$\aleph_1$ such that $\No \M = \aleph_0$'' with \ZFC + \GCH. Namely,
in \cite {Sh125} Shelah proves: it is consistent with \ZFC + \GCH that
there is a group $G$ for which the group of extensions of \Integers by
$G$, in symbols $\Symb{Ext} (G, \Integers)$, is the additive group of
rationals. Here \Integers is the additive group of integers. Then one
extension of \Integers by $G$ can be directly coded to a model \M such
that $\No \M = \Card {\Symb{Ext} (G, \Integers)} = \aleph_0$. The
\Lan-equivalence between two coded models follows from the group
theoretic properties of $G$ \Note {$G$ is strongly $\aleph_1$-free}.
But $\Symb{Ext} (G, \Integers)$ is a divisible group and hence this
coding mechanism is not applicable to the case $1 < \No \M <
\aleph_0$. So there was the problem left if is it consistent to have a
model \M of cardinality $\aleph_1$ for which $1 < \No \M < \aleph_0$.

As Shelah did with the Whitehead problem, we transform \Theorem
{Models} into a question of the nature of pure combinatorial set
theory. The combinatorial problem will be a variant of the
uniformization principles and ladder systems given for example in
\cite {ProperForcing} or \cite {EklofMekler}. As a matter of fact the
more complicated ladder systems used here retrace back to the papers
\cite {Sh98} and \cite {Sh125}.

For the benefit of the reader we sketch the ``standard'' notion of
$(\eta, 2)$-uniformization. For a limit ordinal $\delta < \omega_1$, a
ladder on $\delta$ is a strictly increasing $\omega$-sequence of
ordinals with limit $\delta$. Let $S$ be a set of limit ordinals below
$\omega_1$. A ladder system on $S$ is a function \Function \eta S
{\Functions \omega {\omega_1}} such that each $\eta(\delta)$ is a
ladder on $\delta$. A $2$-colouring on $S$ is a function \Function c S
{\Functions \omega {\Braces {0, 1}}}. For all $\delta \in S$ and $n <
\omega$, a $2$-colouring $c$ on $S$ associates the element $c_{\delta,
n}$ \Note {the $(n +1)$th element of the sequence $c(\delta)$} for
each ``step'' $\eta_{\delta, n}$ of a ladder system $\eta$ on $S$,
hence the name $2$-colouring. A $2$-colouring $c$ on $S$ can be
uniformized if there is a function \Function f {\omega_1} {\Braces {0,
1}} satisfying that for all $\delta \in S$ there is $m < \omega$ such
that for all $n < \omega$, $n > m$ implies $f(\eta_{\delta, n}) =
c_{\delta, n}$. Such a function $f$ is called a uniformizing function
and we say that $c$ is uniform with respect to $\eta$. The $(\eta,
2)$-uniformization holds if every $2$-colouring on $S$ is uniform
w.r.t. $\eta$.

For our purpose we need a different kind of ladder system. The main
difference is that instead of the principle ``all colourings are
uniform'' we want to know what the ``number of nonuniform colourings''
can be. We consider colourings which take values in a field, and hence
we can define a natural equivalence relation for colourings. \Note
{The following definition is from \cite {Sh98}, see also \cite
{EkSh559} where colourings which take values in a group are
considered.} For $2$-colourings $c$ and $d$ on $S$ let $c - d$ be the
$2$-colouring $e$ on $S$ defined for all $\delta \in S$ and $n <
\omega$ by $e_{\delta, n} \in \Braces {0, 1}$ and $(e_{\delta, n} +
d_{\delta, n}) \equiv c_{\delta, n} \pmod 2$. Then $2$-colourings $c$
and $d$ on $S$ are equivalent w.r.t. a ladder system $\eta$ on $S$ if
$c - d$ is uniform w.r.t. $\eta$. The number of pairwise nonequivalent
colourings is the number of equivalence classes of $2$-colourings on
$S$ under the given equivalence relation. But as it is pointed out in
\cite [Theorem 6.2] {Sh98}, for all set $S \Subset \omega_1$ of limit
ordinals and ladder systems on $S$, the number of pairwise
nonequivalent colourings is either 1 or $\geq \AlephExp 0$. In our
transformation of \Theorem {Models} the value of \No \M will
correspond to the number of pairwise nonequivalent colourings. So, all
the cases $1 < \No \M \leq \aleph_0$ are ruled out when only standard
ladder systems are considered.


The main result concerning the combinatorial problem is that for all
finite fields $F$,
 \begin{property}
 it is consistent with \ZFC + \GCH that there are ``good'' ladder
system and ``good'' equivalence for colourings \Note {which take
values in $F$} such that the number of pairwise nonequivalent
colourings is \Card F.
 \end{property}
Recall that all finite fields are of the size $p^m$ with $p$ a prime
number and $m$ a positive integer.

In standard ladders each step is one ordinal. The principal idea of
the ``good'' ladders will be answering to the following simple
question: what happens if each step could be a finite set of ordinals,
or even a ``linear combination'' of standard steps?

In order to make our presentation self contained we give proofs of
some facts which are essentially proved elsewhere \Note {mainly in
\cite {Sh64} and \cite {Sh125}}. In \Subsection {Ladders} we give the
exact definitions for the ``good'' ladder systems, colourings, and
equivalence. In \Subsection {Forcing} we introduce some basic facts
about iterated forcing.

In \Section {Colouring} the combinatorial problem is reformulated in a
precise form and a solution of the problem is presented. Some remarks
concerning generalizations are given in \Subsection {Remarks}. Since
ladder systems and uniformization principles are also used in abelian
group theory and general topology this section may be of independent
interest.

\Section {Models} is devoted to the proof of \Theorem {Models}. We
take a ``good'' ladder system and code each colouring \Col a to a
model \MC a. Then all of the coded models will be \Lan-equivalent, and
moreover, they are isomorphic if and only if the corresponding
colourings are equivalent. So the main result really is a
straightforward consequence of the independence result concerning the
combinatorial problem. The coding technique we have used in the proof
of \Theorem {Models} is a nice trick, and may also be of independent
interest. Hence \Section {Models} is written in a way that if the
reader accepts \Theorem {Colouring} on faith, she or he can read only
\Subsection {Ladders} and then directly proceed to reading \Section
{Models}.


\end{SECTION}

\begin{SECTION} {+} {Preliminaries} {Preliminaries}



For all sets $X, Y, Z$, ordinals $\alpha$ and functions \Function f X
Y:
\begin{itemize}

\item the restriction \Res f Z has the meaning \Res f {(Z \Inter \Dom
f)},

\item \Functions X Y is the set of all functions from $X$ into $Y$,

\item \Functions \alpha Y is the set of all $\alpha$-sequences of
elements in $Y$, and \InitSeq \alpha Y is \BigUnion [\beta < \alpha]
{\Functions \beta Y}.

\end{itemize}
Let $S$ be a subset of a limit ordinal $\mu$ with uncountable
cofinality. The set $S$ is \Def {stationary} in $\mu$ if for all
closed unbounded subsets $C$ of $\mu$, $S \Inter C$ is nonempty. The
set $S$ is \Def {bistationary} in $\mu$ if $S$ is stationary in $\mu$
and $\mu \Minus S$ is also stationary in $\mu$.

\begin{SUBSECTION} {-} {Ladders} {Ladder Systems and Colourings}




Suppose \ModelStructure {F} {+, \Times, 0, 1} is a field. We denote
by \Vec [F] the vector space over $F$ freely generated by \Seq {\Gen
\xi} {\xi < \omega_1}. Suppose $y$ is an element of \Vec [F] and
$e_\xi \in F$ are coefficients such that
 \[
	y = \Sum [\xi < \omega_1] {e_\xi \Gen \xi},
 \]
 where only finitely many of the coefficients are nonzero. \Def
{The support of $y$}, in symbols \Supp y, is the set \Set {\xi <
\omega_1} {e_\xi \not= 0}. For all functions \Function f \mu F such
that $\Supp y \Subset \mu \leq \omega_1$, $f(y)$ is a shorthand for
the following element of $F$,
 \[
	\Sum [\xi < \omega_1] {e_\xi \Times f(\xi)}.
 \]
 A subset $Y$ of \Vec [F] is \Def {unbounded} if for all $\theta <
\omega_1$ there is some $y \in Y$ for which $\theta < \min (\Supp y)$.

\begin{DEFINITION}{Ladders}
 \begin{ITEMS}

\ITEM {ladder}
 A \Vec [F]-ladder on $\delta$, where $\delta < \omega_1$ is a limit
ordinal, is a sequence $\Seq {y_n} {n < \omega}$ of elements in \Vec
[F] such that
 \begin{ITEMS}

	\ITEM {union}
	$\BigUnion [n < \omega] {\Supp {y_n}} \Subset \delta$,

	\ITEM {increasing}
	\Seq [\big] {\min (\Supp {y_n})} {n < \omega} is an
	increasing sequence of ordinals with limit $\delta$, and

	\ITEM {not_cover}
	for all $n < \omega$, $\Supp {y_n} \not\Subset \BigUnion
	[m < n] {\Supp {y_m}}$.

 \end{ITEMS}

\ITEM {system}
 A \Vec [F]-ladder system on $S$, where $S$ is a set of limit ordinals
below $\omega_1$, is a function \LadSys from $S$ into the \Vec
[F]-ladders such that for each $\delta \in S$, $\LadSys (\delta)$ is a
\Vec [F]-ladder on $\delta$.

\ITEM {colouring}
 An $F$-colouring on $S$ is a function from $S$ into \Functions \omega
F. The set of all such colourings is \ColSet [S, F].

 \end{ITEMS}
\end{DEFINITION}

For all $\delta \in S$ and \Vec [F]-ladder systems \LadSys on $S$:
 \begin{itemize}

\item
 the $(n +1)$th element in the $\omega$-sequence $\LadSys (\delta)$ is
denoted by \LadSys [_{\delta, n}];

\item
 \Supp {\LadSys (\delta)} is a shorthand for \BigUnion [n < \omega]
{\Supp {\LadSys [_{\delta, n}]}};

\item
 for a function $f$ with $\Supp {\LadSys (\delta)} \Subset \Dom f$ and
$\Ran f \Subset F$, $f(\LadSys (\delta))$ is a shorthand for the
sequence \Seq [\big] {f(\LadSys [_{\delta, n}])} {n < \omega};

\end{itemize}

When $f$ is a function with $\Dom f = \omega_1$ and $\Ran f \Subset
F$, $f(\LadSys)$ denotes the function from $S$ into \Functions \omega
F which maps each $\delta \in S$ into $f(\LadSys (\delta))$.

\begin{DEFINITION}{Uniform}
 Suppose \LadSys is a \Vec [F]-ladder system on $S$, $\Col a \in
\ColSet [S, F]$, and $D$ is a filter over $\omega$ including all
cofinite subsets of $\omega$, i.e., all subsets $I$ of $\omega$ for
which $\omega \Minus I$ is finite.
 \begin{ITEMS}

\ITEM {almost_all}
 If $\delta \in S$ and $f$ is a function with $\Supp {\LadSys
(\delta)} \Subset \Dom f \Subset \omega_1$ and $\Ran f \Subset F$,
then $f(\LadSys [_{\delta, n}]) = \Col [_{\delta, n}] a$ for almost
all $n < \omega$, or in symbols $f(\LadSys (\delta)) \IsAlmost_D \Col
a (\delta)$, when
 \[
	\Set [\big] {n < \omega}
	{f(\LadSys [_{\delta,n}]) = \Col [_{\delta,n}] a} \in D.
 \]

\ITEM {res}
 If $f$ is a function with $\mu \Subset \Dom f$ and $\Ran f \Subset
F$, then $f$ uniformizes \Res {\Col a} {\mu +1} with respect to
\LadSys and $D$, when $f(\LadSys (\delta)) \IsAlmost_D \Col a
(\delta)$ for all $\delta \in S \Inter \mu +1$.

\ITEM {uniform}
 An $F$-colouring \Col a on $S$ is uniform w.r.t. \LadSys and $D$ if
there is \Function f {\omega_1} F satisfying $f(\LadSys (\delta))
\IsAlmost_D \Col a (\delta)$ for all $\delta \in S$. The set of all
uniform $F$-colourings on $S$ w.r.t. \LadSys and $D$ is \UnifSet
[\LadSys, D].

\ITEM {equivalent}
 The set \ColSet [S, F] forms a vector space over the field $F$, when
addition in \ColSet [S, F] and operation of $F$ on \ColSet [S,F] are
defined componentwise, and the unit element for addition is the
function which is constantly 0. Using the addition of this space we
define \Col a and \Col b in \ColSet [S, F] to be equivalent
w.r.t. \LadSys and $D$, written $\Col a \Equivalent_{\LadSys, D} \Col
b$, if $\Col a - \Col b$ is a uniform colouring w.r.t. \LadSys and
$D$. We denote by \GenBy [F] {\Col a} the subspace of \ColSet [S, F]
generated by $\Col a \in \ColSet [S, F]$.

 \end{ITEMS}
\end{DEFINITION}

It is easy to see that the set \UnifSet [\LadSys, D] forms a subspace
of \ColSet [S, F]. So the factor space \Quotient {\ColSet [S, F]}
{\UnifSet [\LadSys, D]} also forms a vector space over $F$, and
consequently, for all $\Col a, \Col b \in \ColSet [S, F]$, $\Col a
\Equivalent_{\LadSys, D} \Col b$ if and only if \Col a and \Col b
belong to the same coset of \Quotient {\ColSet [S, F]} {\UnifSet
[\LadSys, D]}. If $A$ and $C$ are subsets of \ColSet [S, F] then $A +
C$ is \Set {\Col a + \Col c} {\Col a \in A \And \Col c \in C}. Hence
$\GenBy [F] {\Col b} + \UnifSet [\LadSys, D]$ denotes the set
\ARRAY{
	\Set [\big] {\Col a + \Col c}
	{\Col a \in \GenBy [F] {\Col b} \And \Col c \in
		\UnifSet[\LadSys,D]} \\
	= \Set [\big] {(e \Times \Col b) + \Col c}
	{e \in F \And \Col c \in \UnifSet[\LadSys,D]} \\
	= \Set [\big] {\Col d \in \ColSet [S, F]}
	{\ThereIs e \in F \SuchThat
	e \Times \Col b \Equivalent_{\LadSys, D} \Col d}.
}

\begin{LEMMA} {Colouring}
 Suppose $D$ is a filter over $\omega$ including all cofinite sets of
$\omega$, $S \Subset \omega_1$ is a set of limit ordinals, $F$ is a
field, and \LadSys is a \Vec [F]-ladder system on $S$.
 \begin{ITEMS}

\ITEM {extension}
If \Col a is an $F$-colouring on $S$, $\mu_0 < \omega_1$, and
\Function {f_0} {\mu_0} F uniformizes \Res {\Col a} {\mu_0 +1}
w.r.t. \LadSys and $D$, then for all $\mu_1 < \omega_1 \Minus
(\mu_0+1)$, there is an extension \Function {f_1} {\mu_1} F of $f_0$
which uniformizes \Res {\Col a} {(\mu_1 +1)} w.r.t. \LadSys and $D$.

 \ITEM {nonstationary}
 If $S$ is nonstationary in $\omega_1$, then all $F$-colourings on $S$
are uniform w.r.t. \LadSys and $D$.

 \ITEM {initial}
 Let \Col a be an $F$-colouring on $S$ and $g$ a function from
$\omega_1$ into $F$. If there exists $\mu < \omega_1$ such that
$g(\LadSys (\delta)) \IsAlmost_D \Col a (\delta)$ for all $\delta \in
S \Minus \mu$, then \Col a is uniform w.r.t. \LadSys and $D$.
 \end{ITEMS}



\begin{Proof}%
\ProofOfItem {extension}
 Suppose $S$ is enumerated by \Set {\delta_\alpha} {\alpha <
\omega_1}, where $\delta_\alpha < \delta_\beta$ for all $\alpha <
\beta < \omega_1$, and $e^{\alpha, n}_\xi \in F$ for $\xi, \alpha <
\omega_1$ and $n < \omega$, are coefficients such that
 \[
	\LadSys [_{\delta_\alpha, n}] =
	\Sum [\xi < \delta_\alpha] {e^{\alpha,n}_\xi \Gen \xi}.
 \]

Our first task is to find a function \Function {g_\alpha} {\Supp
{\LadSys (\delta_\alpha)}} F, for all $\alpha < \omega_1$, such that
the equation $g_\alpha (\LadSys [_{\delta_\alpha, n}]) = \Sum [\xi <
\delta_\alpha] e^{\alpha, n}_\xi \Times g_\alpha (\xi) = \Col
[_{\delta_\alpha, n}] a$ holds for all $n < \omega$. Hence consider
the following system of equations,
 \EQUATION{system}{
	 \ForAll n < \omega,
	\Sum [\xi < \delta_\alpha] e^{\alpha, n}_\xi \Times g_\alpha(\xi)
	= \Col [_{\delta_\alpha, n}] a.
 }%
 By \ItemOfDefinition {Ladders} {not_cover} the set $\Supp {\LadSys
[_{\delta_\alpha, n}]} \Minus \BigUnion [m < n] {\Supp {\LadSys
[_{\delta_\alpha, m}]}}$ is nonempty for all $n < \omega$. Besides $F$
is a field. Thus it is possible to define directly by induction on $n
< \omega$ a solution \Function {g_\alpha} {\Supp {\LadSys
(\delta_\alpha)}} F for the system of the equations \Equation
{system}.

We prove by induction on $\alpha < \omega_1$, the following claim,
 \begin{property}
 for all $\mu_0 < \delta_\alpha$ and \Function {f_0} {\mu_0} F
uniformizing \Res {\Col a} {\mu_0 +1}, there is \Function {f_1}
{\delta_\alpha} F uniformizing \Res {\Col a} {\delta_\alpha +1} and
satisfying $f_0 \Subset f_1$.
 \end{property}

Suppose $\mu_0 = 0$ and $\alpha = 0$. Then $f_1 = g_0 \Union \Set
[\big] {(\xi, 0)} {\xi \in \delta_0 \Minus \Dom {g_0}}$ satisfies the
claim.

Suppose $\alpha = \beta +1$, $\mu_0 < \delta_\alpha$, and \Function
{f_0} {\mu_0} F uniformizes \Res {\Col a} {\mu_0 +1}. Let $g_\alpha$
be a solution for the system of the equations \Equation {system}. We
may assume $\mu_0 \geq \delta_\beta$ since if not, then by the
induction hypothesis there is \Function {f_0'} {\delta_\beta} F
extending $f_0$ and uniformizing \Res {\Col a} {\delta_\beta +1}. It
suffices to prove the claim for such $f_0'$.

Define a function \Function {f_1} {\delta_\alpha} F, for all $\xi <
\delta_\alpha$, by
 \EQUATION{successor}{
    \FunctionDefinition {f_1 (\xi)}{
	\FunctionDefMidCase {f_0(\xi)}
			{\If \xi \in \mu_0 = \Dom {f_0}}
	\FunctionDefMidCase {g_\alpha(\xi)}
			{\If \xi \in \Dom {g_\alpha} \Minus \mu_0}
	\FunctionDefOtherwise {0}
    }
 }
 Then of course $f_0 \Subset f_1$ and for all $\delta \in S \Inter
\delta_\alpha = (S \Inter \mu_0) \Union \Braces {\delta_\beta}$, $f_1
(\LadSys (\delta)) = f_0 (\LadSys (\delta)) \IsAlmost \Col a
(\delta)$. By \ItemOfDefinition {Ladders} {increasing} \Set [\big] {n
< \omega} {\Supp {\LadSys [_{\delta_\alpha, n}]} \Inter \delta_\beta
\not= \emptyset} must be finite. Therefore also $f_1 (\LadSys
(\delta_\alpha)) \IsAlmost g_\alpha (\LadSys (\delta_\alpha)) = \Col a
(\delta_\alpha)$ holds. So $f_1$ uniformizes \Res {\Col a}
{\delta_\alpha +1}.

Suppose then $\alpha$ is a limit ordinal. If the limit $\sup (S \Inter
\delta_\alpha) = \theta$ is smaller than $\delta_\alpha$, i.e.,
$\delta_\alpha$ is not a limit of its predecessors in $S$, then we may
assume $\mu_0 = \Dom {f_0} \geq \theta$ by the induction
hypothesis. Furthermore, the function $f_1$ given in \Equation
{successor}, this time for different $\alpha$ of course, is a
uniformizing function for \Res {\Col a} {\delta_\alpha +1}.

Suppose $\delta_\alpha$ is a limit point in $S$, i.e., $\theta =
\delta_\alpha$.  Let \Seq {\epsilon_m} {m < \omega} be an increasing
sequence of ordinals in $S$ with limit $\delta_\alpha$. By the
induction hypothesis there are for all $m < \omega$ functions
\Function {h_m} {\epsilon_m} F uniformizing \Res {\Col a} {\epsilon_m
+1} and satisfying $h_m \Subset h_{m +1}$. This time we may assume
$\Dom {f_0} = \mu_0 = \epsilon_0$ and $f_0 = h_0$. Define a function
\Function {f_1} {\delta_\alpha} F, for all $\xi < \delta_\alpha$, by
 \[
    \FunctionDefinition {f_1 (\xi)}{
	\FunctionDefMidCase {f_0(\xi)}
			{\If \xi < \epsilon_0 = \mu_0 = \Dom{f_0}}
	\FunctionDefMidCase {g_\alpha(\xi)}
			{\If \xi \in \Dom {g_\alpha} \Minus \Dom{f_0}}
	\FunctionDefLastCase {h_l(\xi)}
		{\Text{otherwise, where}
		l = \min \Set {m < \omega} {\xi < \epsilon_m = \Dom {h_m}}}
    }
 \]
 In the definition above, $g_\alpha$ is a solution for \Equation
{system}. Clearly $f_0 \Subset f_1$ and $f_1 (\LadSys (\delta)) = f_0
(\LadSys (\delta)) \IsAlmost \Col a (\delta)$ for all $\delta \in S
\Inter \mu_0$. For all $\delta \in S \Inter \delta_\alpha$, the set
\Set [\big] {n < \omega} {\Supp {\LadSys [_{\delta, n}]} \Inter (\Dom
{f_0} \Union \Dom {g_\alpha}) \not= \emptyset} is finite. Thus for all
$\delta \in S \Inter \delta_\alpha$, there is some $m < \omega$ such
that $f_1 (\LadSys (\delta)) \IsAlmost h_m (\LadSys (\delta))
\IsAlmost \Col a (\delta)$. Since also \Set [\big] {n < \omega} {\Supp
{\LadSys [_{\delta_\alpha, n}]} \Inter \Dom {f_0} \not= \emptyset} is
finite, $f_1 (\LadSys (\delta_\alpha)) \IsAlmost g_\alpha (\LadSys
(\delta_\alpha)) = \Col a (\delta_\alpha)$ holds. So $f_1$ uniformizes
\Res {\Col a} {\delta_\alpha +1}.

\ProofOfItem {nonstationary}
 Suppose \Col a is an $F$-colouring on $S$, and $C = \Set {\mu_\alpha}
{\alpha < \omega_1}$ is a closed and unbounded subset of $\omega_1$
disjoint from $S$. We define by induction on $\alpha < \omega_1$
functions \Function {f_\alpha} {\mu_\alpha} F such that \BigUnion
[\alpha < \omega_1] {f_\alpha} is a uniformizing function for \Col
a. We may assume $\mu_0 = 0$. So let $f_0$ be the function with empty
domain. Suppose $\alpha > 0$ and for all $\gamma < \beta < \alpha$,
functions $f_\gamma$, $f_\beta$, satisfying $f_\gamma \Subset f_\beta$
and $f_\beta$ uniformizing \Res {\Col a} {\mu_\beta +1}, are defined.

If $\alpha$ is a successor of the form $\beta +1$, let \Function
{f_\alpha} {\mu_\alpha} F be some extension of $f_\beta$ which
uniformizes \Res {\Col a} {\mu_\alpha +1}. This is possible by \Item
{extension}. If $\alpha$ is a limit ordinal then $f_\alpha = \BigUnion
[\beta < \alpha] {f_\beta}$ uniformizes \Res {\Col a} {\mu_\alpha +1}
by induction hypothesis, and since $\mu_\alpha \in C \Minus S$. It
follows that $f = \BigUnion [\alpha < \omega_1] {f_\alpha}$
uniformizes \Col a.

\ProofOfItem {initial}
 Suppose \Function g {\omega_1} F satisfies $g(\LadSys (\delta))
\IsAlmost \Col a (\delta)$ for some $\mu < \omega_1$ and for all
$\delta \in S \Minus \mu$. By \Item {extension} there is \Function f
\mu F which uniformizes \Res {\Col a} {\mu +1}. Now, as in the proof
of \Item {extension}, the function $h$ defined for all $\xi <
\omega_1$ by
 \[
    \FunctionDefinition {h(\xi)}{
	\FunctionDefMidCase {f(\xi)} {\If \xi < \mu = \Dom f}
	\FunctionDefMidCase {g(\xi)} {\Otherwise}
 }
 \]
 uniformizes \Col a.
\end{Proof}%
%

\end{LEMMA}

\begin{ShComment}%
\Remark%
\end{ShComment}
\begin{PvComment}%
 A generalization of \Vec [F]-ladder systems and $F$-colourings on $S$
for larger cardinals $\kappa$ is straightforward. Furthermore, a
generalization of \Lemma {Colouring} would hold if $S$ consists of
limit ordinals with fixed cofinality and $S \Inter \mu$ is
nonstationary in $\mu$ for all $\mu < \kappa$.

\end{PvComment}
\begin{ShComment}
 It is possible to replace in \ItemOfDefinition {Ladders} {increasing}
$\min$ by $\max$. It is also possible to replace in \Definition
{Uniform} the filter $D$ by a sequence \Seq {D_\delta} {\delta \in S}
of filters. Such replacements allows more freedom, but in the proof of
\Lemma {Colouring} one should prove by induction the following
slightly stronger statement: if $f_0$ and a finite extension of it
with domain $\subset \mu_1$ are given, then there is an extension
$f_1$ as in \ItemOfLemma {Colouring} {extension}.

On the other hand one may like to replace the field by a ring. In this
case for \Lemma {Colouring} to work it is convenient to demand in
addition to \ItemOfDefinition {Ladders} {ladder} that
 \begin{itemize}

	\item the sets \Supp {y_n}, $n < \omega$, are pairwise
	disjoint, and

	\item for each $n < \omega$, $y_n$ satisfies that for every
	$b$ in the ring $F$ there is a function $f$ with $f(y_n) = b$.

 \end{itemize}
 However, at present work there is no real need for these variants.
\end{ShComment}

\end{SUBSECTION}

\begin{SUBSECTION} {-} {Forcing} {Forcing}




All forcing arguments are considered to be taking place in the
universe $V$ of all sets. Let \ModelStructure P {\Str_P, \One_P} be a
forcing notion, where $\One_P$ is a unique maximal element with
respect to the order $\Str_P$. The subscript $P$ from $\One_P$ will be
omitted everywhere else except in definitions. For all conditions $p$
in $P$, $p \Forces [P] \phi$ means $p$ forces a sentence $\phi$. If
every condition forces $\phi$, we write $\Forces [P] \phi$. The order
$\Str_P$ of conditions $p, q \in P$ is interpreted in a way that $q$
is a stronger condition than $p$ if $q \Str_P p$. Hence for all
sentences $\phi$, $p \Forces [P] \phi$ implies $q \Forces [P] \phi$,
when $q \Str_P p$. The subscript $P$ in the notation $\Str_P$ is not
written when $P$ is obvious from the context.

\begin{SimpleText}%
A subset $X$ of $P$ is an \Def {antichain}, if the elements of $X$ are
pairwise incompatible, i.e., for all $p, q \in X$, there is no $r \in
P$ satisfying both $r \Str_P p$ and $r \Str_P q$. A forcing notion $P$
has the \Def {$\aleph_2$-chain condition}, abbreviated by
$\aleph_2$-c.c., if every antichain of $P$ has cardinality less than
$\aleph_2$. A subset $X$ of $P$ is said to be \Def {dense below $p$},
where $p \in P$, if for all $q \Str_P p$ there is $r \in X$ stronger
than $q$.%
\end{SimpleText}

Let $G$ be a $P$-generic set over $V$. When $\sigma$ is a $P$-name,
the interpretation of $\sigma$ in the generic extension $V [G]$ is
denoted by \Int G \sigma. For an object $o$ in $V[G]$, a $P$-name for
$o$ is written \Name o, i.e., $\Int G {\Name o} = o$. The canonical
name for the generic set $G$ itself is \WideName G. If an object $o$
is in $V$, we identify the name \Name o with the object $o$ itself
instead of using standard names \SimpleCite [page 190] {Kunen}. The
only exceptions for these rules are that the standard names for
uncountable cardinals and collections \Functions Y X are written
\SName [_\alpha] \omega and \WideSName {\Functions Y X} respectively,
to distinguish them from the cardinals $\aleph_\alpha$, $\alpha > 0$,
and corresponding collections in the generic extension. If \Name f is
a $P$-name for a function from $X \in V$ into $Y \in V$ and $x \in X$,
a condition \Def {$p \in P$ decides the value of $\Name f(x)$} when
there is $y \in Y$ satisfying $p \Forces [P] \Name f(x) = y$.

If $P$ is a forcing notion having $\aleph_2$-c.c. then $P$ preserves
all cofinalities $\geq \aleph_2$, i.e., for all limit ordinals
$\theta$, if $\Cf \theta = \kappa \geq \aleph_2$ in $V$ then $\Forces
[P] \Cf \theta = \kappa$. Hence $P$ preserves all cardinals too, i.e.,
if $\lambda \geq \aleph_2$ is a cardinal in $V$ then $\Forces [P]
\Quote {\lambda \Text {is a cardinal}}$.%
 \begin{SimpleText}
 For proofs of these standard facts see \SimpleCite [Lemma 5.10 on
page 206] {Kunen}.
\end{SimpleText}

Suppose that \ModelStructure P {\Str_P, \One_P} is a forcing notion in
$V$ and \WideName Q, \NameStr [_Q], and \Name [_Q] \One are $P$-names
satisfying $\Forces [P] \Quote {\ModelStructure {\WideName Q}
{\NameStr [_Q], \Name [_Q] \One} \Text {is a forcing notion}}$. \Def
{The two stage iteration} \ModelStructure {\CompF P {\WideName Q}}
{\Str_{\CompF P {\WideName Q}}, \One_{\CompF P {\WideName Q}}} is
defined by
 \[
	\CompF P {\WideName Q} =
	\Set [\big] {(p, \Name q)} {p \in P
	\And p \Forces [P] \Name q \in \WideName Q},
 \]
 and for the elements in \CompF P {\WideName Q}, $(p, \Name q)
\Str_{\CompF P {\WideName Q}} (p', \Name q')$ if both $p \Str_P p'$
and $p \Forces [P] (\Name q \NameStr [_Q] \Name q')$ hold. So
$\One_{\CompF P {\WideName Q}}$ is the pair $(\One_P, \One_{\WideName
Q})$. We identify elements $(p, \Name q), (p', \Name q') \in \CompF P
{\WideName Q}$ if both $(p, \Name q) \Str_{\CompF P {\WideName Q}}
(p', \Name q')$ and $(p', \Name q') \Str_{\CompF P {\WideName Q}} (p,
\Name q)$ hold. This iteration amounts to the same generic extension
as does the composition where one first forces with $P$ and then with
\WideName Q. 

\begin{RefArea}{Env}{Definition}{IteratedForcing}

\Def {An iterated forcing of length $\omega_2$ with countable
support},
 \[
	\ModelStructure {P_{\omega_2}}
		{\Str_{P_{\omega_2}}, \One_{P_{\omega_2}}}
	= \CountLim
	\Seq {P_\alpha, \Name [_\alpha] Q} {\alpha < \omega_2}
 \]
 is inductively defined for all $\alpha \leq \omega_2$ as follows.
 \begin{ITEMS}

\ITEM {P_0}
 The forcing notion \ModelStructure {P_0} {\Str_{P_0}, \One_{P_0}} is
defined by $\One_{P_0} = \emptyset$, $P_0 = \Braces {\One_{P_0}}$, and
$\Str_{P_0} = P_0 \times P_0$.

\ITEM {P_alpha}
 Suppose for all $\beta < \alpha$, \Name [_\beta] Q, \NameStr
[_{Q_\beta}], \Name [_{Q_\beta}] \One are given $P_\beta$-names and
they satisfy
\[
	\Forces [P_\beta] \Quote {
		\ModelStructure {\Name [_\beta] Q} {\NameStr [Q_\beta],
		\Name [_{Q_\beta}] \One} \Text{is a forcing notion}
	}.
\]
Moreover, assume that for all $\beta < \alpha$,
\[
	\ModelStructure {P_\beta} {\Str_{P_\beta}, \One_{P_\beta}}
	= \CountLim
		\Seq {P_\gamma, \Name [_\gamma] Q} {\gamma < \beta}
\]
are already defined. It follows from \Item {P_0} that $V = V[H]$ for
all $P_0$-generic sets $H$ over $V$. Hence we assume that \Name [_0]
Q, \NameStr [_{Q_0}], \Name [_{Q_0}] \One are standard names and
\ModelStructure {Q_0} {\Str_{Q_0}, \One_{Q_0}} is a forcing notion in
$V$.

The set $P_\alpha$ is the collection of all functions $p$ satisfying
the following requirements:
 \begin{ITEMS}

\ITEM {domain}
 The domain of $p$ is $\alpha$, and for each $\beta < \alpha$ the
value of $p(\beta)$ is a $P_\beta$-name such that $\Res p \beta
\Forces [P_\beta] p(\beta) \in \Name [_\beta] Q$.

\ITEM {support}
 The set
	\Set {\beta < \alpha}
	{\Res p \beta \DoesNotForce [P_\beta]
	p(\beta) = \Name [_{Q_\beta}] \One}
 is countable.

\end{ITEMS}

\ITEM {Str}
 For all $\alpha \leq \omega_2$ and $p, q \in P_\alpha$, the order of
these conditions is $q \Str_{P_\alpha} p$ if either $\alpha$ is a
limit ordinal, and
 \[
	\ForAll \beta < \alpha, \ 
	\Res q \beta \Str_{P_\beta} \Res p \beta,
 \]
 or otherwise, $\alpha$ is a successor ordinal of the form $\beta +
1$, and
 \ARRAY{
	\Res q \beta \Str_{P_\beta} \Res p \beta, \\
	\Res q \beta \Forces [P_\beta]
	q(\beta) \NameStr [_{Q_\beta}] p(\beta).
 }

\ITEM {One}
 $\One_{P_\alpha}$ is the function which maps each $\beta < \alpha$
into \Name [_{Q_\beta}] \One.

\end{ITEMS}

\end{RefArea} 

\Remark For all $\alpha \leq \omega_2$ and $p \in P_\alpha$, we let
\Dom p denote the set of ordinals given in \DefinitionItem
{IteratedForcing} {support} above. This set is usually called the
support of $p$. So, one can as well think that the domain of a
condition $p \in P_\alpha$ really is the set \Dom p. We may write $f
\in P_\alpha$, $\alpha \leq \omega_2$, when $f$ is only a function
satisfying $\Dom f \Subset \alpha$ and $f \Union \Set {(\beta, \Name
[_{Q_\beta}] \One)} {\beta \in \alpha \Minus \Dom f}$ is a condition
in $P_\alpha$. We abbreviate $\Forces [P_\alpha]$ by $\Forces
[\alpha]$ and $\Str_{P_\alpha}$ by $\Str_\alpha$, or even more
compactly by $\Str$ when the subscript is obvious. 

For each $\beta < \omega_2$, \CompF {P_\beta} {\Name [_\beta] Q} is
isomorphic to $P_{\beta +1}$ via the mapping \Mapping {(p, \Name q)}
{\ConCat p {\SimpleSeq {\Name q}}}. If $G_\alpha$ is a
$P_\alpha$-generic set over $V$ then for each $\beta < \alpha$,
$G_\beta$ denotes the $P_\beta$-generic set \Set {\Res p \beta} {p \in
G_\alpha}.%
 \begin{SimpleText}
 For more details, see for example \SimpleCite [pages 273--275] {Kunen}.
\end{SimpleText}

\begin{FACT} {iteration}
 Suppose $\alpha \leq \omega_2$ and $P_\alpha = \CountLim \Seq
{P_\beta, \Name [_\beta] Q} {\beta <\alpha}$.
 \begin{ITEMS}

\ITEM {aleph2cc}
 If $P_\beta$ has $\aleph_2$-c.c. for all $\beta < \alpha$, then
$P_\alpha$ has $\aleph_2$-c.c.

\ITEM {intermediate}
 If $\alpha = \omega_2$, $P_{\omega_2}$ has $\aleph_2$-c.c., $X$ is a
set in $V$, and \WideName Y is $P_{\omega_2}$-name satisfying $\Forces
[\omega_2] (\WideName Y \Subset X \And \Card {\WideName Y} < \SName
[_2] \omega)$, then for all $P_{\omega_2}$-generic sets $G$ over $V$,
there is $\alpha < \omega_2$ such that the subset $Y = \Int G
{\WideName Y}$ is already in $V[G_\alpha]$.

\ITEM {enumeration}%
Let $S$ be a set of limit ordinals $< \omega_1$ and $F$ a field of
cardinality $\leq \aleph_1$. If $\AlephExp 1 = \aleph_2$ and $\Forces
[\beta] \Par {\Card {\Name [_\beta] Q} = \Card {\SName [_1] \omega}}$
for all $\beta < \alpha$, then there is a collection \Set {\NameCol
[^{\alpha, \gamma}] c} {\gamma < \omega_2} of $P_\alpha$-names
satisfying $\Forces [\alpha] \Set {\NameCol [^{\alpha, \gamma}] c}
{\gamma < \SName [_2] \omega} = \NameColSet [S, F]$. Such a collection
is called $(P_\alpha, \omega_2)$-enumeration for \NameColSet [S, F].

 \end{ITEMS}

\begin{SimpleText}

 \begin{Proof}
 Again these facts are standard forcing arguments. But this time we
give a sketch of proof for \Item {enumeration}.

\ProofOfItem {aleph2cc} \cite [Theorem 2.2] {Baumgartner}

\ProofOfItem {intermediate}
 The iteration $P_{\omega_2}$ preserves the cofinality of $\omega_2$
and thus $\Forces [\omega_2] \Card {\WideName Y} < \Card {\SName [_2]
\omega} = \Cf {\SName [_2] \omega}$. The claim follows from \cite
[Lemma 5.14 on page 276] {Kunen}.

\ProofOfItem {enumeration}
 It suffices to show that $\Card {P_\alpha} \leq \aleph_2$ and
$P_\alpha$ has $\aleph_2$-c.c. \cite [Lemma 5.13 on page 209]
{Kunen}. We prove these properties by induction on $\beta \leq
\alpha$.

If $\beta = 0$ then $P_\beta$ is the trivial forcing. Thus $Q_0$ is in
$V$, $\Card {P_\beta} = \Card {Q_0} = \aleph_1$, and $P_\beta$
trivially satisfies $\aleph_2$-c.c. Assume $0 < \beta \leq \alpha$ and
that the claim holds for all $\gamma < \beta$.

Suppose $\beta$ is a successor ordinal, say $\beta = \gamma +1$. By
induction hypothesis $P_\gamma$ has the $\aleph_2$-c.c. The following
is a well-known fact \cite [Theorem 2.1 and Lemma 3.2] {Baumgartner}:
 \begin{property}
 Suppose $P$ is a forcing notion having $\aleph_2$-c.c. If $\Forces
[P] \Quote {\WideName Q \Text{is a forcing notion with \SName [_2]
\omega-c.c.}}$, then \CompF P {\WideName Q} has $\aleph_2$-c.c.
Furthermore, if $\AlephExp 1 = \aleph_2$, $\Card P \leq \aleph_2$, and
$\Forces [P] \Card {\WideName Q} \leq \SName [_2] \omega$, then $\Card
{\CompF P {\WideName Q}} \leq \aleph_2$.
 \end{property}
 This together with our assumption $\Forces [\gamma] \Par {\Card
{\Name [_\gamma] Q} = \Card {\SName [_1] \omega}}$ and the induction
hypothesis implies \CompF {P_\gamma} {\Name [_\gamma] Q} has
$\aleph_2$-c.c., and moreover $\Card {\CompF {P_\gamma} {\Name
[_\gamma] Q}} \leq \aleph_2$. The claim follows from the knowledge
that $P_\beta$ is isomorphic to \CompF {P_\gamma} {\Name [_\gamma] Q}.

If $\beta$ is a limit ordinal, then $P_\beta = \CountLim \Seq
{P_\gamma, \Name [_\gamma] Q} {\gamma <\beta}$, $P_\gamma$ has the
$\aleph_2$-c.c., and $\Card {P_\gamma} \leq \aleph_2$, for all $\gamma
< \beta$, by the induction hypothesis. Directly by \FactItem
{iteration} {aleph2cc} $P_\beta$ has the $\aleph_2$-c.c. By the
definition, $p \in P_\beta$ iff $\Res p \gamma \in P_\gamma$ for all
$\gamma < \beta$. So the cardinality of the set $P_\beta$ is at most
$\AlephExp [2] 1 = \aleph_2 \Times \AlephExp 1 = \aleph_2$, since
$\Card \beta \leq \aleph_1$ and we assume $\AlephExp 1 = \aleph_2$.
 \end{Proof}

\end{SimpleText}
\end{FACT} 

For $\alpha < \beta \leq \omega_2$, $p \in P_\alpha$ and $q \in
P_\beta$ such that $p \Str_\alpha \Res q \alpha$ the ``composition''
of these conditions, in symbols $\CompCond p q$, is the function
having domain $\alpha$ and defined for all $\gamma < \alpha$ by
 \[
    \FunctionDefinition {(\CompCond p q) (\gamma)}{
	\FunctionDefMidCase {p(\gamma)} {\If \gamma < \beta}
	\FunctionDefLastCase{q(\gamma)} {\If \beta \leq \gamma < \alpha}
 }
 \]
 Then, as in \cite [Definition 1.1 and Fact 1.3] {Sh64} or \cite
[Definition 1.12 and Fact 1.13] {Goldstern}, \CompCond p q is a
condition in $P_\beta$ and $(\CompCond p q) \Str_\beta q$.

We shall also need the ``quotient'' forcing notion \ModelStructure
{\Name [_{\alpha, \beta}] P} {\NameStr [_{\alpha, \beta}], \Name
[_{\alpha, \beta}] \One} of an iterated forcing $P_\beta = \CountLim
\Seq {P_\gamma, \Name [_\gamma] Q} {\gamma < \beta}$, where $\alpha <
\beta \leq \omega_2$. The following definition is from \cite
{Goldstern}. The $P_\alpha$-name \Name [_{\alpha,\beta}] P is such
that
 \[
	\Forces [\alpha] \Name [_{\alpha,\beta}] P =
	\Set {p \in P_\beta} {\Res p \alpha \in \Name [_\alpha] G},
 \] 
 \NameStr [_{\alpha, \beta}] is a $P_\alpha$-name for which
 \[
	\Forces [\alpha]\,
	\NameStr [_{\alpha, \beta}] \  = \,
	\Res {\mathord{\Str_\beta}} {\Name [_{\alpha, \beta}] P},
 \]
 and \Name [_{\alpha, \beta}] \One is the standard name for
$\One_{P_\beta}$. So, for all $P_\alpha$ generic sets $H$ over $V$ and
$p, q \in P_{\alpha, \beta} = \Int H {\Name [_{\alpha, \beta}] P}$, we
have $p \Str_{\alpha, \beta} q$ in $V[H]$ iff $p \Str_\beta q$ in $V$,
where $\Str_{\alpha, \beta}\, = \Int H {\NameStr [_{\alpha,
\beta}]}$. We abbreviate $\Forces [P_{\alpha, \beta}]$ by $\Forces
[\alpha, \beta]$.

\newcommand {\DNNameO} {\Hat{\Hat o}}

\begin{SimpleText}
 We define a mapping \Function {\iota_{\alpha, \beta}} {P_\beta}
{\CompP \alpha \beta} for all $\alpha < \beta \leq \omega_2$, $p \in
P_\beta$ by setting $\iota_{\alpha, \beta} (p) = (\Res p \alpha, p)$.

\begin{FACT} {Names}%
Suppose $\alpha < \beta \leq \omega_2$.
\begin{ITEMS}

\ITEM {embedding}%
$\iota_{\alpha, \beta}$ is a dense embedding.

\ITEM {b->CompP}%
There is a mapping $\iota_{\alpha, \beta}^\star$ from $P_\beta$-names
into \CompP \alpha \beta-names determined by $\iota_{\alpha, \beta}$
such that for all $P_\beta$-names $\Name o$ and \CompP \alpha
\beta-generic sets $G^\star$ over $V$,
 \[
	\Int G {\Name o} =
	\Int {G^\star} {\iota_{\alpha, \beta}^\star(\Name o)},
 \]
 where $G = {\iota_{\alpha, \beta}}^{-1}(G^\star)$. Furthermore, for
all $p \in P_\beta$, formulas $\phi$, and $P_\beta$-names $\Name o$,
 \[
	p \Forces [\beta] \phi(\Name o)
	\Iff
	\iota_{\alpha, \beta} (p) \Forces [\CompP \alpha \beta]
	\phi(\iota_{\alpha, \beta}^\star(\Name o)).
 \]

\ITEM {CompP->ab}%
For each \CompP \alpha \beta-names $\sigma$ there is a $P_\alpha$-name
$\ddot \sigma$ for a \Name [_{\alpha, \beta}] P-name such that for all
$P_\alpha$-generic sets $H$ over $V$ and $P_{\alpha, \beta}$-generic
sets $K$ over $V [H]$,
 \[
	\Int K {\Int H {\ddot \sigma}} = \Int {\CompF H K} \sigma.
 \]

\ITEM {b->ab}%
For each $P_\beta$-names $\Name o$ there is a $P_\alpha$-name
$\DNNameO$ for a \Name [_{\alpha, \beta}] P-name such that for all
$P_\alpha$-generic sets $H$ over $V$ and $P_{\alpha, \beta}$-generic
sets $K$ over $V [H]$,
 \[ 
	\Int K {\NName o} = \Int G {\Name o},
 \]
where $\NName o = \Int H {\DNNameO}$ and $G = {\iota_{\alpha,
\beta}}^{-1} (\CompF H K)$.

 \end{ITEMS}


\end{FACT}
\end{SimpleText}

\begin{FACT} {Quotient}
Suppose $\alpha < \beta \leq \omega_2$, $H$ is a $P_\alpha$-generic
set over $V$, \Name o is a $P_\beta$-name, and $\phi$ is a
formula. Then there is a $P_{\alpha, \beta}$-name \NName o in $V[H]$
such that the following hold.
\begin{ITEMS}

\ITEM {b->ab}
 If $p \in P_\beta$, $\Res p \alpha \in H$, and $p \Forces [\beta]
\phi(\Name o)$ then in $V[H]$, there is $q \in P_{\alpha, \beta}$ such
that $q \Str_{\alpha, \beta} p$ and $q \Forces [\alpha, \beta]
\phi(\NName o)$.

\ITEM {ab->b}
 If in $V[H]$, $r \in P_{\alpha, \beta}$ and $r \Forces [\alpha,
\beta] \phi (\NName o)$ then in $V$, there is $s \in P_\beta$
satisfying $s \Str_\beta r$, $\Res s \alpha \in H$, and $s \Forces
[\beta] \phi(\Name o)$.

 \end{ITEMS}

\begin{SimpleText}
\begin{Proof}%
We shall use the name \NName o given in \ItemOfFact {Names}
{b->ab}. The following is a standard fact \cite [Lemma 23.4 on page
233] {Jech}:
 \TEXTPROPERTY{CompF}{
If $G$ is a $P$-generic set over $V$ and $H$ is a $Q$-generic set in
$V[G]$, then the set
 \[
	\CompF G H =
	\Set [\big] {(p, \Name q) \in \CompF P {\WideName Q}}
	{p \in G \And\, \Int G {\Name q} \in H}
 \]
is a \CompF P {\WideName Q}-generic set over $V$, and $V[G][H] =
V[\CompF G H]$.
}
Another standard property is the following \cite [Theorem 7.11 on
page 221] {Kunen}:
\TEXTPROPERTY{ComP}{%
For all \CompP \alpha \beta-generic sets $G^\star$ over $V$, $G =
{\iota_{\alpha, \beta}}^{-1} (G^\star)$ is a $P_\beta$-generic set
over $V$ and $V[G^\star] = V[G]$.
}

\ProofOfItem {b->ab}%
Take any $P_{\alpha, \beta}$-generic set $K$ over $V[H]$ such that $p
\in K$. Let $G$ be ${\iota_{\alpha, \beta}}^{-1} (\CompF H K)$.  Then
$p \in G$. Moreover, by \Property {CompF}, \Property {ComP}, and
\ItemOfFact {Names} {b->ab},
 \[
	V[H][K] = V[\CompF H K] = V[G]
	\models \phi(o),
 \]
where $o = \Int K {\Int H {\DNNameO}} = \Int {\CompF H K}
{\iota_{\alpha, \beta}^\star (\Name o)} = \Int G {\Name o}$ and
$\DNNameO = (\iota_{\alpha, \beta}^\star (\Name o)) \spddot$. It
follows that there is in $V[H]$ a condition $q' \in P_{\alpha, \beta}$
such that $q' \in K$ and $q' \Forces [\alpha, \beta] \phi(\NName o)$
in $V[H]$. Since $K$ is a filter and $p \in K$ there is $q \in K$ for
which $q \Str_{\alpha, \beta} q', p$. The condition $q$ satisfies the
claim.

\ProofOfItem {ab->b}%
Fix a $P_{\alpha, \beta}$-generic set $K$ over $V[H]$ such that $r
\in K$ in $V[H]$. Then $\Res r \alpha \in H$, $(\Res r \alpha, r) \in
\CompF H K$, and by \Property {CompF} together with \ItemOfFact
{Names} {b->ab},
 \[
	V[H][K] = V[\CompF H K] \models \psi(o),
 \]
where $o = \Int K {\Int H {\DNNameO}} = \Int {\CompF H K}
{\iota_{\alpha, \beta}^\star (\Name o)}$ and $\DNNameO =
(\iota_{\alpha, \beta}^\star (\Name o)) \spddot$. Hence there is
$t^\star \Str_{\CompP \alpha \beta} (\Res r \alpha, r)$ in \CompF H K
such that $t^\star \Forces [\CompP \alpha \beta] \psi (\iota_{\alpha,
\beta}^\star (\Name o))$.

Besides, $\iota_{\alpha, \beta}$ is a dense embedding and hence the
set $\iota_{\alpha, \beta} (P_\beta)$ is dense below $t^\star$.
Consequently, there is $s \in P_\beta$ such that $s \in \CompF H K$
and $\iota_{\alpha, \beta} (s) = (\Res s \alpha, s) \Str_{\CompF H K}
t^\star$. So $\iota_{\alpha, \beta} (s) \Forces [\CompP \alpha \beta]
\psi (\iota^\star (\Name o))$, and by \ItemOfFact {Names} {b->CompP},
$s \Forces [\beta] \psi (\Name o)$. Since $s \Str_\beta r$ and $\Res s
\alpha \in H$, the condition $s$ satisfies the claim.
\end{Proof}

\end{SimpleText}
\end{FACT} 

\begin{FACT}{Res}
 Suppose $\alpha \leq \beta \leq \omega_2$, $p, q \in P_\beta$, and
$H$ is a $P_\alpha$-generic set over $V$. If both $\Res p \alpha \in
H$ and $\Res q \alpha \in H$ hold, then there are $p', q' \in P_\beta$
such that $p' \Str_\beta p$, $q' \Str_\beta q$, and $\EqualRes {p'}
{q'} \alpha \in H$.

\begin{SimpleText}
\begin{Proof}
 Since $H$ is a filter, there is $r \in H$ such that $r \Str_\alpha
\Res p \alpha$ and $r \Str_\alpha \Res q \alpha$. So the conditions
$p' = \CompCond r p$ and $q' = \CompCond r q$ satisfy the claim.
 \end{Proof}

\end{SimpleText}
\end{FACT}


\end{SUBSECTION}


\end{SECTION}

\begin{SECTION} {+} {Colouring} {The Combinatorial Problem}




This section is devoted to the proof of the following theorem which
is a precise form of the theorem described in the introduction.

\begin{THEOREM} {Colouring}
 Assume the following properties hold in $V$:
 \begin{itemize}

\item
 the generalized continuum hypothesis, \GCH;

\item
 $S$ is a set of limit ordinals below $\omega_1$ and bistationary in
$\omega_1$;

\item
 $F$ is a finite field;

\item
 \Vec is the vector space over $F$ freely generated by \Seq {\Gen \xi}
{\xi <\omega_1};

\item
 $D$ is a filter over $\omega$ including all cofinite sets of $\omega$.

 \end{itemize}
 Then there is a forcing notion \ModelStructure {P} {\Str, \One} of
cardinality $\aleph_2$ such that $P$ satisfies $\aleph_2$-c.c., $P$
does not add new countable sequences, and for every $P$-generic set
$G$ over $V$, there is in $V[G]$ a \Vec-ladder system \LadSys on $S$
such that $\Card {\Quotient {\ColSet [S, F]} {\UnifSet [\LadSys, D]}}
= \Card F$.
 \end{THEOREM}

Recall that the conclusion of the theorem is equivalent to the number
of pairwise nonequivalent $F$-colourings on $S$ w.r.t. \LadSys and $D$
being \Card F. The idea of the forthcoming proof of the theorem will
be similar to the proof of \cite [Theorem 1] {Sh125}.

\relax From now on, all \Vec-ladders on $\delta$ and \Vec-ladder systems on
$S$ are called simply \Def {ladders on $\delta$} and \Def {ladder
systems}, all $F$-colourings on $S$ are called \Def {colourings} for
short, and \ColSet denotes the set of all $F$-colourings on $S$. The
subspace of \ColSet generated by a colouring \Col b is shortly \GenBy
{\Col b}.

\begin{SUBSECTION} {-} {Definition} {Definition of the Forcing}




To define an iterated forcing $P = \CountLim \Seq {P_\alpha, \Name
[_\alpha] Q} {\alpha < \omega_2}$ it suffices to define names for
forcing notions \ModelStructure {\Name [_\alpha] Q} {\NameStr
[_{Q_\alpha}], \Name [_{Q_\alpha}] \One} by induction on $\alpha <
\omega_2$.

The forcing notion \ModelStructure {Q_0} {\Str_{Q_0}, \One_{Q_0}} is
defined as follows. The set $Q_0$ is $\InitLadSet \times \InitColSet$
where
 \ARRAY[lll]{
	\InitLadSet &=&
	\Set {\Res {\LadSysSymbol z} \theta}
	{\LadSysSymbol z \Text {is a ladder system and}
	\theta < \omega_1}, \\
	\InitColSet &=&
	\Set {\Res {\Col c} \mu} {\Col c \in \ColSet \And \mu < \omega_1}.
 }
 We shorten our notation for $p = (\Res {\LadSysSymbol z} \theta, \Res
{\Col c} \mu) \in Q_0$ by writing
 \ARRAY[l]{

\text{p[1] for \Res {\LadSysSymbol z} \theta
and p[2] for \Res {\Col c} \mu,}\\

\text{$\epsilon \leq \Dom p$ if $\epsilon \leq \Min {\theta, \mu}$, and}\\

\text{$\Dom p \leq \epsilon$ if $\Max {\theta, \mu} \leq \epsilon$.}

 }
 For all $p_0, p_1 \in Q_0$, we define $p_1 \Str_{Q_0} p_0$ iff $p_1$
coordinatewise extends $p_0$, i.e., $p_1[1] \Superset p_0[1]$ and
$p_1[2] \Superset p_0[2]$. The pair of functions with empty domain is
the maximal element $\One_{Q_0}$ of $Q_0$. If $X \Subset Q_0$ is a set
of pairwise compatible conditions then we define
 \[
	\BigCompCond {\Set p {p \in X}} =
	\big(\BigUnion {\Set {p[1]} {p \in X}},
	 \BigUnion {\Set {p[2]} {p \in X}}\big).
 \]
 Note that $Q_0$ is $\aleph_1$-closed \Note {which means every
descending $\omega$-chain of conditions has a lower bound}. Hence
$Q_0$ does not add new countable sequences and $\aleph_1$ is not
collapsed\SimpleCite [Theorem 6.14 on page 214] {Kunen}.

For every $P_1$-generic set $G_1$ there are $P_\alpha$-names \Name
\LadSys and \NameCol b, for $\alpha = 1$ \Note {later on $\alpha$
might be any index in $\omega_2 \Minus \Braces 0$}, such that
 \ARRAY{
 \Forces [\alpha]
   \Name \LadSys =
	\BigUnion {\Set {p(0)[1]} {p \in \Name [_\alpha] G}} \\

 \Forces [\alpha]
  \NameCol b = \BigUnion {\Set {p(0)[2]} {p \in \Name [_\alpha] G}}.

 }
 So, these names together with a generic set determine a ladder system
and a colouring. Hereafter \Def {uniform} and \Def {equivalent} mean
uniform and equivalent w.r.t. the generic ladder system \Name \LadSys
and the filter $D$. Hence \UnifSet denotes the set of all uniform
colourings w.r.t. \Name \LadSys and $D$. Observe that the generic
colouring \NameCol b satisfies $\Forces [1] (\NameCol b \not\in
\NameUnifSet)$, as we shall prove in \Lemma {Nonuniform}.

Forcing notions \ModelStructure {\Name [_\alpha] Q} {\NameStr
[_{Q_\alpha}], \Name [_{Q_\alpha}] \One}, for $1 \leq \alpha <
\omega_2$, are defined in such a way that each \Name [_\alpha] Q
``kills'' an undesirable colouring. In order to ensure that all
undesirable colourings will be killed, a bookkeeping function will be
needed. Fix $\pi$ to be a function from $\omega_2$ onto $\omega_2
\times \omega_2$ such that whenever $\pi(\alpha) = (\beta, \gamma)$
then $\beta \leq \alpha$.

The bookkeeping function is useful only if we can ensure that the
colourings can be enumerated by $\omega_2$. Since we assume \GCH the
cardinality of \ColSet is $\Card {\Functions S {(\Functions \omega
F)}} = (\AlephExp 0)^{\aleph_1} = \AlephExp 1 = \aleph_2$. Hence there
is an enumeration \Set {\Col [^{0, \gamma}] c} {\gamma < \omega_2} for
\ColSet in $V$. By \ItemOfFact {iteration} {enumeration} the existence
of a $(P_\alpha, \omega_2)$-enumeration for \NameColSet follows for $1
\leq \alpha < \omega_2$, if we show that for each $\beta < \alpha$,
 \PROPERTY{cardQ}{
	\Forces [\beta]
	\Card {\Name [_\beta] Q} \leq \Card {\SName [_1] \omega}.
 }%
 Since $P_0$ is the trivial forcing \Braces \One, $\AlephExp 0 =
\aleph_1$, and $\Card \InitColSet = \Card \InitLadSet = (\AlephExp
0)^{\aleph_0} = \AlephExp [1] 0$, we have that $\Card {Q_0} =
\aleph_1$, and so \Property {cardQ} holds trivially when $\beta = 0$.

Suppose $1 \leq \alpha < \omega_2$. Our induction hypothesis is that
for each $\beta < \alpha$, there is a $(P_\beta,
\omega_2)$-enumeration \Set {\NameCol [^{\beta, \gamma}] c} {\gamma <
\omega_2} for \NameColSet and that $\Forces [\beta] (\Card {\Name
[_\beta] Q} = \Card {\SName [_1] \omega})$ holds. It follows from
\ItemOfFact {iteration} {enumeration} that there also exists a
$(P_\alpha, \omega_2)$-enumeration \Set {\NameCol [^{\alpha, \gamma}]
c} {\gamma < \omega_2} for \NameColSet.

\begin{DEFINITION} {Q_alpha}
 Suppose $\pi(\alpha) = (\beta, \gamma)$. Then $\beta \leq \alpha$ and
\NameCol [^{\beta,\gamma}] c has been defined. We define \NameCol
[^\alpha] a to be a $P_\alpha$-name which refers to the same colouring
as the $P_\beta$-name \NameCol [^{\beta, \gamma}] c, i.e., for every
$P_\alpha$-generic sets $H$ over $V$, $\Int H {\NameCol [^\alpha] a} =
\Int {H_\beta} {\NameCol [^{\beta, \gamma}] c}$. A $P_\alpha$-name
\Name [_\alpha] Q is defined by
 \[
 \Forces [\alpha]
   \FunctionDefinition {\Name [_\alpha] Q}{
	\FunctionDefMidCase {\Braces {\Name [_{Q_\alpha}] \One}}
		{\If \NameCol [^\alpha] a \in
		\GenBy {\NameCol b} + \NameUnifSet}
	\FunctionDefMidCase {\NameUnifFunc [^\alpha] a} {\Otherwise}
   }
 \]
 where \Name [_{Q_\alpha}] \One is the standard name for the function
having empty domain, and \NameUnifFunc [^\alpha] a is a $P_\alpha$-name
satisfying
 \[
	\Forces [\alpha] \NameUnifFunc [^\alpha] a =
	\Set {f} {\mu < \SName [_1] \omega \And \Function f \mu F
	\Text{uniformizes} \Res {\NameCol [^\alpha] a} {\mu+1}}.
 \]
 A $P_\alpha$-name \NameStr [_{Q_\alpha}] is defined by $\Forces
[\alpha] (\ForAll p, q \in \Name [_\alpha] Q,\ p \NameStr
[_{Q_\alpha}] q \Iff p \Superset q)$.
 \end{DEFINITION} 

For every $p \in P_\alpha$, an index $\beta \leq \alpha$ is called
\Def {$p$-trivial} if $\beta > 0$ and $\Res p \beta \Forces [\beta]
\Name [_\beta] Q = \Braces \One$. Observe that if $\beta \in \Dom p$
then $\Res p \beta \DoesNotForce [\beta] \Par {p(\beta) = \One}$, and
$\beta$ is not $p$-trivial. Note also that $\Forces [\alpha] \Par
{\NameUnifFunc [^\alpha] a \not= \Braces \One}$ by \ItemOfLemma
{Colouring} {extension}. In fact, if $p \in P_\alpha$ and $p$ forces
\Par {\NameCol [^\alpha] a \not\in \GenBy {\NameCol b} + \NameUnifSet}
then $p$ forces \Name [_\alpha] Q to be a nontrivial forcing notion
\Note {see \ItemOfLemma {Basic} {killing} below}.

We have to check that the property \Property {cardQ} for $\beta =
\alpha$ holds. We shall prove that $P_\alpha$ does not add new
countable sequences. Hence $\Forces [\alpha] \WideSName {\InitSeq
{\SName [_1] \omega} F} = \InitSeq {\SName [_1] \omega} F$. This
implies that
 \[
	\Forces [\alpha]
	\Card {\Name [_\alpha] Q} \leq
	\Card [\big] {\InitSeq {\SName [_1] \omega} F} =
	\Card [\big] {\WideSName {\InitSeq {\SName [_1] \omega} F}} =
	\Card {\SName [_1] \omega},
 \]
 since $\Card {\InitSeq {\omega_1} F} = \AlephExp 0 = \aleph_1$.

Before proving that $P_\alpha$ does not add new countable sequences,
we introduce useful notations and lemmas. Let \HModel \beta {}, for
$\beta \leq \alpha$, denote the model
 \[
	\ModelStructure [\Big] {H(\lambda)}{
	\in, \beta, S, F, D,
	\Seq [\big] {\ModelStructure {P_\gamma} {\Str, \One}}
		{\gamma \leq \beta}
	},
 \]
 where $\lambda$ is ``some large enough'' cardinal, for example
${(\beth_{\omega_2})}^+$, and $H(\lambda)$ is the set of all sets
hereditary of cardinality $< \lambda$. The expansion of the model
\HModel \beta {} with new constant symbols ``$X_1, X_2, \dots$'' is
denoted by \HModel \beta {X_1, X_2, \dots}.

A condition $p$ in $P_\beta$ has \Def {height $\epsilon$}, where
$\beta \leq \alpha$ and $\epsilon < \omega_1$, if for every $\gamma
\in \Dom p$, $\Res p \gamma \Forces [\gamma] \Dom {p(\gamma)} =
\epsilon$. We say that \Def {$p$ is of height $< \epsilon$} when $\Res
p \gamma \Forces [\gamma] \Dom {p(\gamma)} < \epsilon$. The notion $p$
is of height $\geq \epsilon$ is defined analogously. These notions are
from \cite {Sh125}.

If $X$ is a set of pairwise compatible conditions in $P_\alpha$, the
``composition'' of these conditions, in symbols \BigCompCond [(p \in
X)] p, is the function $f$ with $\Dom f = \BigUnion [p \in X] {\Dom
p}$ and for each $\beta \in \Dom f$, $f(\beta)$ is a $P_\beta$-name
such that
 \[
 \Forces [\beta]
  \FunctionDefinition {f(\beta)}{
	\FunctionDefMidCase {\BigCompCond {\Set {p(0)} {p \in X}}}
		{\If \beta = 0}
	\FunctionDefOtherwise {\BigUnion {\Set {p(\beta)} {p \in X}}}
 }
 \]
 Observe that $f$ is not necessarily a condition in $P_\alpha$
\Note{as we pointed out earlier, by this we mean that not even the
extended function $f \Union \Set {(\beta, \One)} {\beta \in \alpha
\Minus \Dom f}$ is a condition in $P_\alpha$}.

\begin{LEMMA} {UnionAndHeight}
 \begin{ITEMS}

\ITEM {union}
 Suppose $\beta \leq \alpha$, \Seq {p_n} {n < \omega} is a descending
chain of conditions in $P_\beta$, $\theta < \omega_1$ is a limit
ordinal not in $S$, and \Seq {\theta_n} {n < \omega} is an increasing
sequence of ordinals with limit $\theta$. Suppose also that for all
$\gamma < \beta$,
\begin{ITEMS}

\ITEM {geq}
 there are infinitely many $m < \omega$ for which $\Res {p_m} \gamma
\Forces [\gamma] \Dom {p_m(\gamma)} \geq \theta_m$, and

\ITEM {leq}
 there are infinitely many $n < \omega$ such that $\Res {p_n} \gamma
\Forces [\gamma] \Dom {p_n(\gamma)} \leq \theta$.
 
\end{ITEMS}
Then $q = \BigCompCond [n < \omega] {p_n}$ is a condition in
$P_\beta$, $q \Str p_n$ for every $n < \omega$, and $q$ has height
$\theta$.

\ITEM {height}
 For all $\beta \leq \alpha$, $p \in P_\beta$, and $\epsilon <
\omega_1$ there are $q \Str p$ in $P_\beta$ and $\theta < \omega_1$
such that $\epsilon \leq \theta$ and $q$ has height $\theta$.

 \end{ITEMS}



\begin{Proof}%
The idea of the proof is similar to \cite [Lemma 1.5] {Sh64}.

\ProofOfItem {union}
 We prove the claim by induction on $\beta \leq \alpha$. If $\beta =1$
then $q \in P_1 \in V$, and clearly the other properties hold
too. Suppose $\beta > 1$ and for every $\gamma < \beta$, $\Res q
\gamma \in P_\gamma$, $\Res q \gamma \Str \Res {p_n} \gamma$ for all
$n < \omega$, and \Res q \gamma has height $\theta$. If $\beta$ is a
limit ordinal then the claim holds directly by the definition of
$P_\beta$ and height. Note that \Dom q is countable even if $\beta$
has cofinality $> \omega$ since \Dom q is a countable union of
countable sets.

Suppose $\beta = \gamma +1$ and $\gamma \in \Dom q$ \Note {if $\gamma
\not\in \Dom q$ then the claim follows from the induction hypothesis}.
By the definition of $q$, $\Res q \gamma \Forces [\gamma] \BigUnion [n
< \omega] {p_n(\gamma)} = q(\gamma)$. By \Item {leq} and \Item {geq},
$\Res q \gamma$ forces that $\Dom {q(\gamma)} = \BigUnion [m < \omega]
{\theta_m} = \theta$. Since $\theta \not\in S$ and $\Res q \gamma
\Forces [\gamma] p_n(\gamma) \in \Name [_\gamma] Q$,
 \[
	\Res q \gamma \Forces [\gamma] \Quote {
	\BigUnion [n < \omega] {p_n(\gamma)} = q(\gamma)
	\Text {uniformizes} \Res {\NameCol [^\gamma] a} {\theta+1}
	}.
 \]
 Consequently, $q \in P_\beta$, $q \Str p_n$ for all $n < \omega$, and
$q$ has height $\theta$.

\ProofOfItem {height}
 Again we work by induction on $\beta \leq \alpha$. If $p \in P_1$ and
$0 \in \Dom p$ then any extension $q \in P_1$ of $p$ for which $\Dom
{q(0)} \geq \epsilon$ suffices to prove the claim. Suppose $\beta =
\gamma +1$, $\gamma \in \Dom p$, and as the induction hypothesis, $r
\Str_\gamma \Res p \gamma$ is a condition in $P_\gamma$ having height
$\theta (\geq \epsilon)$. Since $\Res p \gamma \Wkr_\gamma r \Forces
[\gamma] \Par {p(\gamma) \in \Name [_\gamma] Q}$ we get by
\ItemOfLemma {Colouring} {extension} that $r$ forces
 \[
	\ThereIs x \in \Name [_\gamma] Q \ForWhich
	x \NameStr [_{Q_\gamma}] p(\gamma) \And \Dom x \geq \theta.
 \]
 By the Maximal Principle \SimpleCite [Theorem 8.2 on page 226]
{Kunen} there is a $P_\gamma$-name \Name f satisfying the formula
above and moreover, we may assume $r \Forces [\gamma] \Dom {\Name f} =
\theta$. Define a condition $q \in P_\beta$ by $\Res q \gamma = r$ and
$q(\gamma) = \Name f$. Then $q$ has height $\theta$.

Suppose that $\beta$ is a limit ordinal, and for all $p' \in P_\beta$,
$\gamma < \beta$, and $\epsilon' < \omega_1$ there is a condition $r$
in $P_\gamma$ satisfying $r \Str \Res {p'} \gamma$ and $r$ has height
$\theta' \geq \epsilon'$. We assume that the supremum of \Dom p is
$\beta$ \Note {otherwise the claim follows by the induction
hypothesis}. We define by induction on $n < \omega$ a descending chain
\Seq {q_n} {n < \omega} of conditions in $P_\beta$ such that $q =
\BigCompCond [n < \omega] {q_n}$ will be a condition in $P_\beta$ and
$q$ has height $\theta (\geq \epsilon)$.

Let \Seq {\gamma_n} {n < \omega} be an increasing sequence of ordinals
with limit $\beta$ \Note {$\beta = \sup (\Dom p)$ must be of
cofinality $\omega$}. Note that the set of all $\theta < \omega_1$,
for which
 \begin{property}
	there is a countable elementary submodel
	\EM of \HModel [n < \omega] \beta {p, \gamma_n}
	such that $\EM \Inter \omega_1 = \theta$,
 \end{property}
 is closed and unbounded in $\omega_1$. Because $S$ is bistationary in
$\omega_1$ we can choose a countable elementary submodel \EM of the
model \HModel [n < \omega] \beta {p, \gamma_n} for which $\EM \Inter
\omega_1 = \theta \geq \epsilon$ and $\theta \not\in S$. Let \Seq
{\epsilon_n} {n < \omega} be an increasing sequence of ordinals with
limit $\theta$ \Note {$\epsilon_n \in \EM$ for every $n <
\omega$}. The model \EM satisfies our induction hypothesis and $p,
\gamma_0 \in \EM$, thus there is a condition $r_0 \Str \Res p
{\gamma_0}$ in $P_{\gamma_0} \Inter \EM$ having height greater than
$\epsilon_0$. We define $q_0$ to be $\CompCond {r_0} p$ \Note {which
really is a condition in $P_\beta \Inter \EM$}. Similarly, when the
condition $q_n \in P_\beta \Inter \EM$ is defined we can find a
condition $q_{n +1} \in P_\beta \Inter \EM$ such that $q_{n +1}
\Str_\beta q_n$ and the initial segment \Res {q_{n +1}} {\gamma_{n
+1}} has height greater than $\epsilon_{n +1}$. So \Item {geq} holds
for \Seq {q_n} {n < \omega} and \Seq {\epsilon_n} {n < \omega}. Since
the conditions $q_n$, $n < \omega$, are in \EM and $\EM \Inter
\omega_1 = \theta$, also \Item {leq} is satisfied. It follows from
\Item {union} that $q = \BigCompCond [n < \omega] {q_n}$ is a
condition in $P_\beta$ having height $\theta (\geq \epsilon)$.
\end{Proof}%
%

\end{LEMMA}

Now we are ready to show that $P_\alpha$ is $\aleph_1$-distributive
\Note {see the next lemma}. Hence it will follow that $\aleph_1$ is
not collapsed and for every $P_\alpha$-generic sets $G_\alpha$ over
$V$, if $X \in V$ and $V[G_\alpha] \models (\Function f \mu X \And \mu
< \omega_1)$, then $f$ is already in $V$.

\begin {LEMMA} {Distributive}
 If $E_n$, $n < \omega$, are dense and open subsets of $P_\alpha$,
then \BigInter [n < \omega] {E_n} is dense.

\begin{Proof}%
 \begin{SimpleText}%
 First of all, $E_n$ dense means that for each $p \in P_\alpha$ there
is $q \in E_n$ with $q \Str p$, and $E_n$ open means that if $p \in
E_n$ and $q \Str p$ then $q \in E_n$. Suppose $p$ is an arbitrary
condition in $P_\alpha$. We show that there is a condition $q$ in
\BigInter [n < \omega] {E_n} stronger than $p$.

\end{SimpleText}
 Let \EM be a countable elementary submodel of \HModel [n < \omega]
\alpha {p, E_n} for which $\EM \Inter \omega_1 = \epsilon \in
\omega_1$ and $\epsilon \not\in S$ \Note {for the existence of such
model, see the proof of \ItemOfLemma {UnionAndHeight} {height}}. Fix
an increasing sequence \Seq {\epsilon_n} {n < \omega} of ordinals with
limit $\epsilon$. We define by induction on $n < \omega$ conditions
$q_n \in P_\alpha$ such that for each $n < \omega$,
\ARRAY{
	q_n \in E_n, \\
	q_n \Text {is of height} \geq \epsilon_n, \\
	q_n \Wkr q_{n+1}.
}
Since \EM is an elementary submodel, $E_0 \Inter \EM$ is a dense
subset of $P_\alpha \Inter \EM$. So there is a condition $r \in E_0
\Inter \EM$ stronger than $p$. We let $q_0$ be some extension of $r$
having a height greater than $\epsilon_0$. This is possible since
$\epsilon_0$ is in \EM, and \EM is an elementary submodel of \HModel
[n < \omega] \alpha {p, E_n} which satisfies \ItemOfLemma
{UnionAndHeight} {height}. Moreover, $q_0$ is in $E_0$ since $E_0
\Inter \EM$ is an open subset of $P_\alpha \Inter \EM$. Similarly, if
$q_n \in P_\alpha \Inter \EM$ is already defined we can find $q_{n+1}
\in P_\alpha \Inter \EM$ satisfying the properties given above.

As in the proof of \ItemOfLemma {UnionAndHeight} {height}, $q =
\BigCompCond [n < \omega] {q_n}$ really is a condition in
$P_\alpha$. Now $q \Str q_n$ for each $n < \omega$, and since $E_n$,
$n < \omega$, are open sets, it follows that $q \in \BigInter [n <
\omega] {E_n}$.
 \end{Proof}
\end{LEMMA}

\relax From the preceding lemma it follows that for all $\alpha \leq
\omega_2$ and $p \in P_\alpha$ there is $q \Str p$ in $P_\alpha$
satisfying the following property: for every $\beta < \alpha$, \Res q
\beta decides the value of $q(\beta)$ \Note {proof of this fact can be
made using the same kind of induction as the proof of \ItemOfLemma
{UnionAndHeight} {height}}. Hence, from now on, the reader can think,
if he or she wants, that all conditions in $P_\alpha$ are ``real''
functions from $\alpha$ into \InitSeq {\omega_1} F, not only
``normal'' conditions with names for sequences. Especially, this
thought might me helpful during the first reading of \Lemma {Bound}
below. But we shall use the following conventions. We write $\Dom
{p(\beta)} = \epsilon$, where $p \in P_\alpha$, $\alpha \leq
\omega_2$, $\beta \in \Dom p \Minus \Braces 0$, and $\epsilon \in
\omega_1$, when $p$ is a condition which satisfies $\Res p \beta
\Forces [\beta] \Dom p(\beta) = \epsilon$.  Similarly, we write $\xi
\in \Dom {p(\beta)}$ if $\Res p \beta \Forces [\beta] \Par {\xi \in
\Dom {p(\beta)}}$, and for $c \in F$ we write $p (\beta) (\xi) = c$ if
$\xi \in \Dom {p(\beta)}$ and $\Res p \beta \Forces [\beta]
p(\beta)(\xi) = c$.

We define $g_\alpha$, for nonzero $\alpha < \omega_2$, to be the
generic function determined by \Name [_\alpha] Q, i.e., \Name
[_\alpha] g is a $P_{\alpha +1}$-name satisfying
 \[
	\Forces [\alpha+1] {\Name [_\alpha] g}
	= \BigUnion { \Set {p(\alpha)} {p \in \WideName G} }.
 \]
  Then $g_\alpha$ is a function in $V[H]$ for any $P_{\alpha
+1}$-generic set $H$ since $H$ contains only compatible
conditions. Note that in $V[H]$, $g_\alpha$ is the function with empty
domain iff $Q_\alpha \not= \UnifFunc [^\alpha] a$.

\begin {LEMMA} {Basic}
 \begin{ITEMS}

\ITEM {aleph2cc}
 The forcing notion $P$ is of cardinality $\aleph_2$, and it satisfies
$\aleph_2$-c.c.

\ITEM {distributive}
 $P$ does not add new countable sequences.

\ITEM {collapse}
 For every $P$-generic set $G$ over $V$, $V[G]$ satisfies \GCH and
\Par {(\aleph_\alpha)^V = \aleph_\alpha} for all ordinals $\alpha$.

\ITEM {killing}
 For all nonzero $\alpha < \omega_2$ and $P_{\alpha+1}$-generic sets
$G_{\alpha +1}$ over $V$, $V[G_{\alpha +1}] \models \Col [^\alpha] a
\in \GenBy {\Col b} + \UnifSet$.

\ITEM {bookkeeping}
 For every $P$-generic set $G$ over $V$, $V[G] \models \Card
{\Quotient \ColSet \UnifSet} \leq \Card F$.

 \end{ITEMS}

\begin{Proof}
 Even though all the properties are standard we sketch proofs for
them.

\ProofOfItem {aleph2cc}
 The claim follows directly by the property \PropertyOnPage [Down]
{cardQ} and \ItemOfFact {iteration} {aleph2cc}.

\ProofOfItem {distributive}
 If we assume that there is a new subset of $\omega$ in $V[G]$, where
$G$ is a $P$-generic set over $V$, then by the $\aleph_2$-c.c.
property of $P$ and \ItemOfFact {iteration} {intermediate} we can
choose $\alpha < \omega_2$ such that the new subset is already in
$V[G_\alpha]$. This contradicts \Lemma {Distributive}.

\ProofOfItem {collapse}
 The generalized continuum hypothesis is preserved by \Item
{aleph2cc}, \Item {distributive}, and by the following well-known fact
\SimpleCite [a generalization of Lemma 5.13 on page 209] {Kunen}:
 \begin{property}
 if $\Card P \leq \aleph_2$, $P$ has $\aleph_2$-c.c., $\AlephExp 1 =
\aleph_2$, $\lambda$ is an uncountable cardinal, and $\theta =
({\aleph_2}^\lambda)^V$, then $\Forces [P] 2^\lambda \leq \theta$.
 \end{property}

By \Item {aleph2cc} the ordinals ${\aleph_\alpha}^V$, $\alpha \geq 2$,
are cardinals in the generic extension. Since by \Item {distributive},
${\aleph_1}^V$ is not collapsed, the claim follows.

\ProofOfItem {killing}
 Let $G_{\alpha +1}$ be a $P_{\alpha +1}$-generic set over $V$. If
$(Q_\alpha = \Braces \One)$ holds in $V[G_\alpha]$ then by \Definition
{Q_alpha} $V[G_\alpha] \models \Col [^\alpha] a \in \GenBy {\Col b} +
\UnifSet$. Since $V[G_{\alpha +1}] \Superset V[G_\alpha]$, the latter
formula is also satisfied in $V[G_{\alpha +1}]$.

Suppose \Par {Q_\alpha = \UnifFunc [^\alpha] a} holds in
$V[G_\alpha]$. By \ItemOfLemma {Colouring} {extension} for each $\xi <
\omega_1$ the generic set $G_{\alpha +1}$ contains a condition $p$ for
which $\Res p \alpha \Forces [\alpha] \xi \in \Dom {p(\alpha)}$. Thus
$\Dom {g_\alpha} = \omega_1$ in $V[G_{\alpha +1}]$. Let $f_p$ be a
shorthand for \Int {G_{\alpha +1}} {p (\alpha)}. Then $f_p$
uniformizes \Res {\Col [^\alpha] a} {(\Dom {f_p} +1)} in $V[G_{\alpha
+1}]$. Consequently, $g_\alpha = \BigUnion {\Set {f_p} {p \in
G_\alpha}}$ uniformizes \Col [^\alpha] a in $V[G_{\alpha +1}]$. So
$V[G_{\alpha +1}] \models \Col [^\alpha] a \in \GenBy {\Col b} +
\UnifSet$.

\ProofOfItem {bookkeeping}
 Assume the claim fails. Since $\Card {\GenBy [F] {\Col b}} \leq \Card
F$, let $G$ be a $P$-generic set over $V$ and \Col d a colouring in
$V[G]$ for which $\Col d \not\in \GenBy {\Col b} + \UnifSet$. Since
$P$ has $\aleph_2$-c.c. and $\Forces [P] \Par {\Card {\Col d} < \SName
[_2] \omega}$ there must be, by \ItemOfFact {iteration}
{intermediate}, $\beta < \omega_2$ such that $\Col d \in V
[G_\beta]$. By the definition of the forcing $P$ and \ItemOfFact
{iteration} {enumeration}, \Par {\Set {\Col [^{\beta, \gamma}] c}
{\gamma < \omega_2} = \ColSet} holds in $V[G_\beta]$. So there is
$\gamma < \omega_2$ with $V[G_\beta] \models \Col d = \Col [^{\beta,
\gamma}] c$. By \Definition {Q_alpha} and since the bookkeeping
function $\pi$ is surjective, there is $\alpha < \omega_2$ such that
$(\Col [^\alpha] a = \Col [^{\beta, \gamma}] c)$ holds in
$V[G_\alpha]$. Then by \Item {killing}, $V[G_{\alpha +1}]$ satisfies
$\Col [^\alpha] a \in \GenBy {\Col b} + \UnifSet$. Since $V[G_\alpha
+1] \Subset V[G]$, \Par {\Col [^{\beta, \gamma}] c = \Col [^\alpha] a
= \Col d \in \GenBy {\Col b} + \UnifSet} holds in $V[G]$ contrary to
our initial assumption.
 \end{Proof}
\end{LEMMA}

\Remark It can be seen from the constructions in \Subsection
{Nonuniform} below that $P$ is a proper forcing notion \cite [Theorem
2.8(1) on page 86] {ProperForcing}. But this fact does not, however,
help with the main problem of \Subsection {Nonuniform}.


\end{SUBSECTION}


\begin{SUBSECTION} {+} {Nonuniform}
		{The Generic Colouring is Nonuniform}




The main problem left after \Lemma {Basic} is that maybe the size of
\Quotient \ColSet \UnifSet is smaller than the size of $F$ in the
generic extension. Since $\Card {\Quotient \ColSet \UnifSet} < \Card
F$ implies $\ColSet = \UnifSet$, we may, equivalently, suspect that
the generic colouring \NameCol b is uniform in the generic extension.
As a preliminary lemma we want to show that the generic colouring
\NameCol b is initially nonuniform, but first we have to prove the
following auxiliary lemma.

\begin{LEMMA} {Extension}
 \begin{ITEMS}

\ITEM {P_1}
 Suppose $p \in P_\alpha$, $\alpha \leq \omega_2$, $\delta \in S$, and
$\Dom {p(0)} \leq \delta$. If \Ladder y is a ladder on $\delta$, and
$\bar c$ is an $\omega$-sequence of elements in $F$, then there is $q
\Str p$ satisfying
 \ARRAY{
	\Dom q = \Dom p \Union \Braces 0, \\
	\EqualRes p q {(\alpha \Minus \Braces 0)}, \\
	q \Forces [\alpha] \Name \LadSys (\delta) = \Ladder y
	\And \NameCol b (\delta) = \bar c.
 }

\ITEM {P_alpha}
 Suppose $p \in P_\alpha$, $\alpha \leq \omega_2$, $A$ is a finite
subset of $\alpha \Minus \Braces 0$, \Seq {c_\beta} {\beta \in A} is a
sequence of elements in $F$, and \Seq {y_\beta} {\beta \in A} is a
sequence of elements in \Vec such that $\Supp {y_\beta} \not\Subset
\Dom {p(\beta)}$. Then there is a condition $s \Str p$ in $P_\alpha$
satisfying for all $\beta \in A$ that
 \begin{property}
	either $\beta$ is $s$-trivial or
	$s (\beta) (y_\beta) = c_\beta$.
 \end{property}
 Furthermore, if for each $\beta \in A$,
 \PROPERTY {P_alpha} {
	\Res p \beta \Forces [\beta]
	\Name [_\beta] Q = \NameUnifFunc [^\beta] a,
 }%
 then we can also ensure that
 \ARRAY{
	\Dom s = \Dom p \Union A, \\
	\EqualRes p s {(\alpha \Minus A)}, \\
	\Dom {s (\beta)} = \max(\Supp {y_\beta}) +1.
 }

 \end{ITEMS}

\begin{Proof} This proof is essentially the same as the proof of \cite
[Lemma 1.5] {Sh64}.

 \ProofOfItem {P_1}%
Define $r \in Q_0$ to be any extension of $p(0)$ which satisfies $\Pro
{r[1]} \delta = \Ladder y$ and $\Pro {r[2]} \delta = \bar c$. Then $q$
defined by $\Dom q = \Braces 0$ and $q(0) = r$ is a condition in
$P_1$. Moreover, $q \Str_1 \Res p 1$ and thus the condition \CompCond
q p is as required in the lemma.

\ProofOfItem {P_alpha}
 It suffices to prove the lemma when $A$ is a singleton \Braces \beta,
since the result for larger sets follows by induction, of course
different induction depending on \Property {P_alpha}.

If \Property {P_alpha} holds then define $q = p$, otherwise let $q
\Str p$ in $P_\alpha$ be such that either $\beta$ is $q$-trivial or
\Res q \beta forces \Name [_\beta] Q to be nontrivial. If $\beta$ is
$q$-trivial then $s = q$ is as wanted. Otherwise, assume \Res q \beta
forces \Name [_\beta] Q to be nontrivial. Let $\theta$ be $\max(\Supp
{y_\beta})$. By \ItemOfLemma {Colouring} {extension} \Note {as in the
proof of \ItemOfLemma {UnionAndHeight} {height}} there is a
$P_\beta$-name \Name f for which
 \[
	\Res q \beta \Forces [\beta]
	\Name f \in \Name [_\beta] Q,
	q(\beta) \Subset \Name f \And \theta \Subset \Dom {\Name f}.
 \]
 Define \Name g to be a $P_\beta$-name for a function such that $\Res
q \beta \Forces [\beta] \Par {\Dom {\Name g} = \theta +1, \EqualRes
{\Name f} {\Name g} \theta, \And \Name g (y_\beta) = c_\beta}$. Then
 \[
	\Res q \beta \Forces [\beta] \Name g \Text {uniformizes}
	\Res {\NameCol [^\beta] a} {\theta+2}.
 \]
 Thus \Res q \beta forces both $(\Name g \in \Name [_\beta] Q)$ and
\Par {\Name g \NameStr [_{Q_\beta}] q(\beta)}, and we can define a
condition $r \in P_{\beta +1}$ by $\Dom r = (\Dom q \Inter \beta)
\Union \Braces \beta$, \EqualRes q r \beta, and $r(\beta) = \Name
g$. Then $r \Str_{\beta +1} \Res q {\beta +1}$, and hence $s =
\CompCond r q$ is a condition in $P_\alpha$ satisfying the properties
required.
 \end{Proof}
\end{LEMMA}

\begin{LEMMA} {Nonuniform}
 The generic colouring \NameCol b satisfies $\Forces [1] \NameCol b
\not\in \NameUnifSet$.

\begin{Proof}
 Suppose, contrary to the claim, that there is a condition $p \in P_1$
and $P_1$-name \Name h for a function from $\omega_1$ into $F$ such
that $p$ forces $\Name h (\Name \LadSys) \Equivalent \NameCol b$. Let
\EM be a countable elementary submodel of $\HModel 1 {} (p, \Name h)$
such that $\EM \Inter \omega_1$ is an ordinal $\delta \in S$ \Note
{such \EM exists by a same kind of argument as in the proof of
\ItemOfLemma {UnionAndHeight} {height}}. Choose two increasing
sequences \Seq {\epsilon_n} {n < \omega} and \Seq {\xi_n} {n < \omega}
of ordinals with limit $\delta$. We define by induction on $n <
\omega$ conditions $q_n \in P_1 \Inter \EM$ and elements $d_n \in F$
\Note {$F = F \Inter \EM$ since $F$ is finite}.

Let $r \in P_1 \Inter \EM$ be such that $r \Str p$ and $\epsilon_0
\leq \Dom r$. We define $q_0 \in P_1 \Inter \EM$ to be an extension of
$r$ which decides the value of $\Name h(\xi_0)$, say $d_0 \in F$ and
$q_0 \Forces [1] \Name h(\xi_0) = d_0$. Similarly, if we assume that
$q_n$ and $d_n$ are already defined, we let $q_{n+1} \in P_1 \Inter
\EM$ and $d_{n+1}$ be such that $\epsilon_{n +1} \leq \Dom {q_{n +1}}$
and $q_{n+1} \Forces [1] \Name h(\xi_{n +1}) = d_{n +1}$.

Since $q_n \in \EM$, $\Dom {q_n(0)} < \delta$ holds for every $n <
\omega$. As pointed out many times before, $q = \BigCompCond [n <
\omega] {q_n}$ is a condition in $P_1$ which does not yet decide the
values of $\Name \LadSys (\delta)$ or $\NameCol b (\delta)$. These
properties together with \ItemOfLemma {Extension} {P_1} and the fact
that \Seq {\Gen {\xi_n}} {n < \omega} is a ladder on $\delta$ ensure
that there is $r \Str q$ in $P_1$ satisfying for each $n < \omega$
that $r \Forces [1] \Name [_{\delta,n}] \LadSys = \Gen {\xi_n} \And
\NameCol [_{\delta, n}] b = d_n + 1$. This contradicts the fact that
$r \Str p$ and $p \Forces [1] \Par {\Name h(\Name \LadSys (\delta))
\IsAlmost \NameCol b (\delta)}$, since for all $n < \omega$,
 \[
	r \Forces [1]
	\Name h(\Name [_{\delta,n}] \LadSys)
	= \Name h(\xi_n) = d_n
	\not= d_n +1 = \NameCol [_{\delta,n}] b.
 \]
 \end{Proof} 
\end{LEMMA}

Note that it follows from \ItemOfLemma {Colouring} {nonstationary} and
\Lemma {Nonuniform} that after forcing with the first step $P_1$ the
set $S$ is still stationary in $\omega_1$. An analogous situation also
concerns the forthcoming proof of the theorem: we shall show that
\NameCol b is nonuniform after forcing with the whole iteration $P$,
thus the set $S$ must remain stationary in $\omega_1$ \Note {recall
that cardinals are preserved by \ItemOfLemma {Basic} {collapse}}.

To prove the theorem it suffices to show that the following holds,
 \[
	\Forces [P] \Quote {\NameCol b \Text {is nonuniform}}.
 \]
 Assume, contrary to this claim, that there exists a $P$-generic set
$G$ over $V$ and in the generic extension $V[G]$ a uniformizing
function \Function h {\omega_1} F for the colouring $\Col b = \Int G
{\NameCol b}$. Since $\Card h < \aleph_2$ we can choose, by
\ItemOfLemma {Basic} {aleph2cc} and \ItemOfFact {iteration}
{intermediate}, the minimal ordinal $\alpha^* < \omega_2$ such that
$h$ is already in $V[G_{\alpha^*}]$ \Note {$\alpha^* \geq 2$ by \Lemma
{Nonuniform}}. For the rest of this section, i.e., for the rest of the
proof of \Theorem {Colouring}, let \Name h be a $P_{\alpha^*}$-name,
and $p^* \in P_{\alpha^*}$ be a condition such that
\PROPERTY{contra}{%
	p^* \Forces [\alpha^*]%
	\Quote {\Name h \Text {uniformizes} \NameCol b}.%
}%
By assuming this we are aiming at a contradiction. Note that $G$ is
not fixed. To shorten our notation, we abbreviate the set \Set {p \in
P_{\alpha^*}} {p \Str_{\alpha^*} p^*} by $P^*$. Purely for technical
reasons we assume $0 \in \Dom {p^*}$.

Although the proof of \Lemma {Nonuniform} was simple, it has already
revealed the main idea of the forthcoming proof. Namely, we want to
contradict \Property {contra} by finding an index $\delta^* \in S$ and
a condition $r$ in $P^*$ which forces $\Name h(\Name \LadSys
(\delta^*)) \not\IsAlmost \NameCol b (\delta^*)$. The next lemma
indicates that this is not a trivial task.

\begin{LEMMA} {not_decide}
 If $Y$ is an unbounded subset of \Vec and $d$ is an element in $F$,
then there is no single condition $p \in P^*$ which forces \Par {\Name
h(y) \not= d} for every $y \in Y$.

\begin{Proof}%
Assume such an unbounded set $Y$ and a condition $p \in P^*$
exist. Let \EM be a countable elementary submodel of \HModel
{\alpha^*} {p, Y} such that $\EM \Inter {\omega_1} = \delta \in
S$. Since \EM is an elementary submodel, $Y \Inter \EM$ must be
unbounded in $\delta$. Fix a ladder \Seq {y_n} {n < \omega} on
$\delta$ such that $y_n \in Y \Inter \EM$ for all $n < \omega$. Since
$p \in \EM$ and $\EM \Inter \omega_1 = \delta$, $\Dom {p(0)} <
\delta$. By \ItemOfLemma {Extension} {P_1} there is $q \Str p$ in
$P^*$ satisfying for all $n < \omega$,
\[
	q \Forces [\alpha^*] \Name [_{\delta,n}] \LadSys = y_n 
	\And \NameCol [_{\delta,n}] b = d.
\] 
Since $q \Str p$, $q$ forces \Par {\Name h(\Name [_{\delta, n}]
\LadSys) \not= \NameCol [_{\delta, n}] b}, for all $n < \omega$.  This
contradicts $q \Str p^*$ and $p^*$ forces $\Name h(\Name \LadSys
(\delta)) \IsAlmost \NameCol b (\delta)$.
 \end{Proof}
\end{LEMMA}

Because there is no single condition which decides enough about $\Name
h$ we shall use a descending chain \Seq {p_n} {n < \omega} of
conditions and a lower bound $r$ of the chain. Since $P_\alpha$, for
$2 \leq \alpha \leq \alpha^*$, are not $\aleph_1$-closed, it is not
easy to find suitable chain and bound. The following lemma, together
with \Lemma {TreeSystemBound} and \Lemma {TreeSystem}, solves this
problem. The idea behind the following \LemmaNo {Bound}, \DefinitionNo
{Am-tree}, \DefinitionNo {TreeSystem}, and \LemmaNo {TreeSystemBound}
is similar to the constructions in the proof of \cite [Theorem 1.1]
{Sh64}.



Before the lemmas we fix some notation. Suppose a function $f$ is
\BigCompCond [k < \omega] {p_k} where \Seq {p_k} {k < \omega} is a
descending chain of conditions in $P^*$. Such a function $f$ is said to
be a \Def {countable union of conditions in $P^*$}, and as in \Lemma
{UnionAndHeight}, $f$ has \Def {height $\epsilon$}, where $\epsilon <
\omega_1$, if
\begin{itemize}

\item for each $k < \omega$, $p_k$ is of height $< \epsilon$, and

\item for all $\alpha \in \Dom f$ and $\theta < \epsilon$, there is
$k < \omega$ such that $\alpha \in \Dom {p_k}$ and $p_k$ is of height
$\geq \theta$.
\end{itemize}
For all $\alpha < \alpha^*$, $\xi < \omega_1$, and $c \in F$, we
write $f (\alpha) (\xi) = c$, when there is $n < \omega$ such that
$p_n (\alpha) (\xi) = c$.
 So if $\Ladder y = \Seq {y_n} {n < \omega}$ is a sequence of elements
in \Vec, $\bar a = \Seq {a_n} {n < \omega}$ is a sequence of elements
in $F$, and $\alpha \in \Dom f$ then $f (\alpha) (\Ladder y) \IsAlmost
\bar a$ means that
\[
	\Set [\big] {n < \omega}
	{\ThereIs k < \omega \SuchThat
	\Res {p_k} \alpha \Forces [\alpha]
	p_k(\alpha)(y_n) = a_n} \in D.
\]
We write $f \Subset p$, where $p \in P_\alpha$ and $\alpha \leq
\alpha^*$, if $\Dom f \Subset \Dom p$ and for each $\beta \in \Dom f$
the condition \Res p \beta forces $f(\beta) \Subset p(\beta)$. Note
that if $\alpha \in \Dom f$ then there is $n < \omega$ such that
$\alpha \in \Dom {p_n}$ and $\Res {p_n} \alpha \DoesNotForce [\alpha]
p_n (\alpha) = \One$. It follows that \Res {p_n} \alpha forces \Name
[_\alpha] Q to be nontrivial, and hence $\alpha$ is not $p_m$-trivial
for any $m < \omega$.

Let $\delta^*$ be an ordinal satisfying $\Dom {p^*(0)} < \delta^* \in
S$ and $A^*$ a nonempty and countable subset of $\alpha^* \Minus
\Braces 0$. Suppose $\Braces 0 \Union A^*$ is enumerated by \Set
{\alpha_i} {i < i^*}, where $2 \leq i^* < \omega_1$ and $0 = \alpha_0
< \alpha_i < \alpha_j$ for all $0 < i < j < i^*$.

\begin{LEMMA} {Bound}
 Suppose that $\Ladder y = \Seq {y_n} {n < \omega}$ is a ladder on
$\delta^*$ and for each \Function u {i^*} {\Functions \omega F} there
exists a mapping $f_u$ satisfying the following properties:
 \begin{ITEMS}

\ITEM {union}
 $f_u$ is a countable union of conditions in $P^*$, $\Dom {f_u}
\Subset \Braces 0 \Union A^*$, and $f_u$ has height $\delta^*$;

\ITEM {tree}
 for all \Function {u, v} {i^*} {\Functions \omega F} and $i < i^*$,
if \EqualRes u v i then $\Res {f_u} {\alpha_i} = \Res {f_v}
{\alpha_i}$;

\ITEM {almost}
 for every nonzero $i < i^*$, if $\alpha_i \in \Dom {f_u}$ then $f_u
(\alpha_i) (\Ladder y) \IsAlmost u(i)$.

 \end{ITEMS}
 Then there is \Function u {i^*} {\Functions \omega F} and a condition
$r \in P^*$ such that $f_u \Subset r$, i.e., $r$ is a lower bound for
the conditions which form $f_u$. Moreover, $r$ forces \Par {\Name
\LadSys (\delta^*) = \Ladder y} and $(\NameCol [_{\delta^*, n}] b
\not= 0)$ for every $n < \omega$.




\begin{Proof}%
The proof below is directly based on \cite [Lemma 1.7] {Sh64}.

First of all we define for each \Function u {i^*} {\Functions \omega
F} a condition $r^u_0 \in P_1$ as follows. By \Item {union} $f_u$ is a
union of conditions and $\Dom {f_u(0)} = \delta^*$. Hence, by the
definition of $Q_0$, $\Res {f_u} {\alpha_1} = \Res {f_u} 1$ is a
condition in $P_1$ \Note {$\Dom {\Res {f_u} {\alpha_1}} = \Braces
{\alpha_0} = \Braces 0$}. By \ItemOfLemma {Extension} {P_1} there is a
condition $r^u_0 \Str_1 \Res {f_u} 1$ in $P_1$ for which
 \PROPERTY{0}{
	r^u_0 \Forces [1]
	\Name \LadSys (\delta^*) = \Ladder y \And
	\NameCol [_{\delta^*,n}] b = 1, \ForAll n < \omega.
 }%
 Since $f_u$ is a union of conditions stronger than $p^*$, $r^u_0
\Str_1 \Res {p^*} 1$. Clearly, $\Res {f_u} {\alpha_1} \Subset
r^u_0$. Note that for all \Function {u, v} {i^*} {\Functions \omega F}
if \EqualRes u v 1 then \EqualRes {f_u} {f_v} 1, by \Item {tree}.
Hence we may assume $r^u_0 = r^v_0$ for all $u, v$ satisfying
\EqualRes u v 1.

For technical reasons we define $\alpha_{i^*}$ to be $\alpha_{(i^*
-1)} +1$ if $i^*$ is a successor ordinal and $\sup \Set {\alpha_i} {i
< i^*}$ otherwise. We prove by induction on $k \leq i^*$ the following
extension property for all $1 \leq j < k \leq i^*$:
 \begin{property} 

if \Function u {i^*} {\Functions \omega F} and $p \in P_{\alpha_j}$
satisfy
 \[
	\Res p 1 \Str_1 r^u_0 \And \Res {f_u} {\alpha_j} \Subset p,
 \]
 then there are \Function v {i^*} {\Functions \omega F} and $r \in
P_{\alpha_k}$ such that
 \[
	\EqualRes u v j,
	\Res r {\alpha_j} \Str_{\alpha_j} p,
	\And \Res {f_v} {\alpha_k} \Subset r.
 \]
 \end{property} 

Suppose first that $1 \leq j < k \leq i^*$, $k$ is a successor
ordinal, and \Function u {i^*} {\Functions \omega F} and $p \in
P_{\alpha_j}$ are as required above. Observe that this includes the
case $j = 1$ and $k = j +1 = 2$. We may assume $k = j +1$ since
otherwise there are, by the induction hypothesis, $u'$ extending $u$
and $p'$ such that \EqualRes u {u'} j and $\Res {f_{u'}} {\alpha_{k
-1}} \Subset p'$. It suffices to prove the claim for such $u'$ and
$p'$.

If $\alpha_j \not\in \Dom {f_u}$ then $v = u$ and $r = p$ satisfy the
claim. Assume $\alpha_j \in \Dom {f_u}$. Let $q \Str p$ in
$P_{\alpha_j}$ and a sequence $\bar d \in \Functions \omega F$ be such
that
 \PROPERTY{decides}{
	q \Forces[\alpha_j]
	\NameCol [^{\alpha_j}] a (\delta^*) = \bar d.
 }%
 Note that by \Lemma {Distributive}, $\bar d$ is in $V$. Define a
function \Function v {i^*} {\Functions \omega F} for all $i < i^*$ by
 \[
    \FunctionDefinition {v(i)}{
	\FunctionDefMidCase {\bar d} {\If i = j}
	\FunctionDefOtherwise {u(i)}
    }
 \]
 Since \EqualRes v u j, it follows from \Item {tree} that $\EqualRes
{f_v} {f_u} {\alpha_j} \Subset p \Wkr q$. Let \Seq {p_m} {m < \omega}
be a descending chain of conditions exemplifying that $f_u$ is union
of conditions in $P^*$ and $f_u$ has height $\delta^*$. Then $\Res
{p_m} {\alpha_j} \Wkr_{\alpha_j} q$ for every $m < \omega$, and
furthermore, for each $\delta \in S \Inter \delta^*$ there is $m <
\omega$ such that
 \[
	\Res {p_m} {\alpha_j} \Forces [\alpha_j]
	  f_v (\alpha_j) (\Name \LadSys (\delta))
	= p_m (\alpha_j) (\Name \LadSys (\delta))
	\IsAlmost \NameCol [^{\alpha_j}] a (\delta).
 \]
 By \Item {almost} and since $q \Str_{\alpha_j} \Res {p_m} {\alpha_j}$
the set \Set {n < \omega} {f_v (\alpha_j) (y_n) = v(j)(n)} is in $D$.
This together with $\Res q 1 \Str_1 r^u_0$, \Property {0}, and
\Property {decides} imply that
 \[
	q \Forces [\alpha_j]
	  f_v (\alpha_j) (\Name \LadSys (\delta^*))
	= f_v (\alpha_j) (\Ladder y)
	\IsAlmost v(j) = \bar d
	= \NameCol [^{\alpha_j}] a (\delta^*).
 \]
 We define $r$ to be $q \Union \Braces {(\alpha_j, f_v (\alpha_j))}$.
Then $r$ is a condition in $P_{\alpha_k}$ satisfying $\Res r
{\alpha_j} = q \Str_{\alpha_j} p$ and $\Res {f_v} {\alpha_k} \Subset
r$.

The second case is that $k \leq i^*$ is a limit ordinal. Suppose $1
\leq j < k$ and $u$, $p$ satisfy the assumptions of the extension
property. Our induction hypothesis is that the extension property
holds for all $k' < k$. Let \EM be a countable elementary submodel of
 \[
	\HModel {\alpha^*} {} \Big(
	\delta^*, i^*,
	\Seq {\alpha_i} {i < i^*},
	p, u,
	\Seq {r^w_0} {\Function w {i^*} {\Functions \omega F}},
	\Seq {f_w} {\Function w {i^*} {\Functions \omega F}}
	\Big),
 \]
 such that $\EM \Inter \omega_1 = \theta \in \omega_1 \Minus S$. We
let \Seq {\theta_n} {n < \omega} be an increasing sequence of ordinals
with limit $\theta$, and $\Seq {j_n} {n < \omega}$ be an increasing
sequence of ordinals with limit $k$, where $j_0 = j$. Note that each
$j_n$ is in \EM since $i^* < \omega_1$ and $\M \Inter \omega_1$ is an
ordinal.

We define by induction on $n < \omega$ conditions $q_n \in
P_{\alpha_{j_n}} \Inter \EM$ and functions \Function {u_n} {i^*}
{\Functions \omega F} in \EM as follows. Let $u_0$ be $u$ and $q_0 \in
P_{\alpha_{j_0}} \Inter \EM$ be an extension of $p$ having height
greater than $\theta_0$. This is possible by \ItemOfLemma
{UnionAndHeight} {height}.

Suppose $u_n \in \EM$ and $q_n \in P_{\alpha_{j_n}} \Inter \EM$ are
already defined. Suppose also that $q_n$ has height greater than
$\theta_n$, $\Res {q_n} 1 \Str_1 r^{u_n}_0$, $\Res {f_{u_n}}
{\alpha_{j_n}} \Subset q_n$, and \EqualRes {u_n} {u_m} {j_m} for every
$m < n$. Since \EM is an elementary submodel, our induction hypothesis
holds in \EM. Hence there are in \EM a function $u_{n +1}$ and $r'$ in
$P_{\alpha_{j_{n+1}}}$ with \EqualRes {u_{n +1}} {u_n} {j_n}, $\Res
{r'} {\alpha_{j_n}} \Str_{\alpha_{j_n}} q_n$, and $\Res {f_{u_{n +1}}}
{\alpha_{j_{n+1}}} \Subset r'$. We define $q_{n +1}$ in
$P_{\alpha_{j_{n+1}}} \Inter \EM$ to be an extension of $r'$ having
height greater than $\theta_{n +1}$. Again, this is possible by
\ItemOfLemma {UnionAndHeight} {height}.

Now $\Res {q_{n +1}} {\alpha_{j_n}} \Str_{\alpha_{j_n}} q_n$ and
\EqualRes {u_{n +1}} {u_n} {j_n} for all $n < \omega$. We define $r$
to be \BigCompCond [n < \omega] {q_n}. This is a condition in
$P_{\alpha_k}$ by \ItemOfLemma {UnionAndHeight} {union}. We define a
function \Function v {i^*} {\Functions \omega F} for all $i < i^*$ by
 \[
    \FunctionDefinition {v(i)}{
	\FunctionDefMidCase {u_m(i)} {\If i < k,
		\Where m = \min \Set {n <\omega} {i < j_n}}
	\FunctionDefOtherwise {u(i)}
    }
 \]
 Then directly by their definition and \Item {tree}, $r$ and $v$
satisfy
 \[
	\Res {f_v} {\alpha_k} =
	\BigCompCond [n < \omega] {\Res {f_{u_n}} {\alpha_{j_n}}}
	\Subset \BigCompCond [n < \omega] {q_n} = r.
 \]
\end{Proof}%
%

\end{LEMMA}

Consequently, there is a lower bound for a certain descending chain of
conditions if the functions $f_u$, \Function u {i^*} {\Functions
\omega F}, satisfying the requirements of the preceding lemma exist
\Note {remember, $f_u$ is a union of conditions but not necessarily a
condition itself}. We shall find those functions as unions of
conditions in special kinds of trees. We again need some more
notation. Let $\Filtration A = \Seq {A_m} {m < \omega}$ be a chain of
finite subsets of the set $A^*$ such that $A_m = A^*$ for all $m <
\omega$ if $A^*$ is finite, and otherwise \Filtration A is increasing
and $A^* = \BigUnion [m < \omega] {A_m}$. Such a chain \Filtration A
is called a \Def {filtration of $A^*$}. The disjoint union \BigUnion
[l \leqslant m] {A_l \times \Braces l}, for $m < \omega$, is
abbreviated by \DU A m. For $m < \omega$, $\DU A m \Inter \alpha$ is a
shorthand for the set $\BigUnion [l \leqslant m] {(A_l \Inter \alpha)
\times \Braces l}$, and for a function $\eta$ having the domain \DU A
m, \Res \eta \alpha is a shorthand for the restriction \Res \eta {(\DU
A m \Inter \alpha)}.

\begin{DEFINITION} {Am-tree}
 Suppose $m < \omega$. We set
 \[
	\Ind A m =
	\Set {\eta} {\eta \Text{ is a function from \DU A m into} F}.
 \]
 An \DU A m-condition tree $T$ is a mapping from \Ind A m into $P^*$
with the property that for all $\eta, \nu \in \Ind A m$ and $\alpha
\in A_m$,
 \[
	\EqualRes \eta \nu \alpha \Implies
	\EqualRes {T(\eta)} {T(\nu)} \alpha.
 \]
 Sometimes we abbreviate $T(\eta)$ by $T_\eta$.

Suppose $n \leq m < \omega$. An \DU A m-condition tree $T$ is stronger
than an \DU A n-condition tree $R$, in symbols $T \Str R$, if for each
$\eta \in \Ind A m$, $T (\eta) \Str_{\alpha^*} R (\Res \eta {\DU A
n})$.

An \DU A m-condition tree $T$ is of height $\geq \epsilon$, $\epsilon
< \omega_1$, if all the conditions in $T$ are of height $\geq
\epsilon$. The notion ``\,$T$ has height $< \epsilon$'' is defined
analogously.
 \end{DEFINITION}

\begin{DEFINITION} {TreeSystem}
 Suppose \Filtration A is a filtration of $A^*$, \Ladder y is a ladder
on $\delta^*$, and $\bar \epsilon$ is an increasing sequence of
ordinals with limit $\delta^*$. An $(\bar \epsilon, \Ladder y)$-tree
system on \Filtration A is a family $\TreeSystem T = \Seq {T^m} {m <
\omega}$ of functions fulfilling the following requirements for each
$m < \omega$:
 \begin{ITEMS}

\ITEM {Tree}
 $T^m$ is an \DU A m-condition tree;

\ITEM {domain}
 for all $\eta \in \Ind A m$, $\Dom {T^m_\eta} \Subset \Braces 0
\Union A^*$ \Note {where $A^* = \BigUnion [m < \omega] {A_m}$};

\ITEM {m_value}
 for all $\eta \in \Ind A m$ and $\alpha \in A_m$,
	$\alpha$ is $T^m_\eta$-trivial or
	$T^m_\eta (\alpha) (y_m) = \eta (\alpha, m)$;

\ITEM {height}
 $T^m$ is of height $\geq \epsilon_m$ and $< \delta^* \Note {= \sup
\bar \epsilon$};

\ITEM {stronger}
 $T^m \geq T^{m +1}$.

 \end{ITEMS}
\end{DEFINITION}

Recall that we assume $\alpha \in \Dom {T^m_\eta}$ and $\Res
{T^m_\eta} \alpha \Forces [\alpha] \Supp {y_m} \Subset \Dom {T^m_\eta
(\alpha)}$ when we write $T^m_\eta (\alpha) (y_m) = \eta (\alpha, m)$.

\begin{LEMMA} {TreeSystemBound}
 For each $(\bar \epsilon, \Ladder y)$-tree system \TreeSystem T on
\Filtration A there are indices $\eta^m \in \Ind A m$, $m < \omega$,
such that \Seq {T^m (\eta^m)} {m < \omega} is a descending chain of
conditions having a lower bound $r \in P^*$. Moreover, $r$ forces \Par
{\Name \LadSys (\delta^*) = \Ladder y} and for all $n < \omega$,
$\NameCol [_{\delta^*, n}] b \not= 0$.



\begin{Proof}%
The idea of the following proof is similar to \cite [Lemma 1.8]
{Sh64}. Recall that \Set {\alpha_i} {i < i^*} is an increasing
enumeration of $\Braces 0 \Union A^*$.

For all $m < \omega$ and \Function u {i^*} {\Functions \omega F} we
define the index $\eta^m_u \in \Ind A m$ by setting for all $(\alpha,
n) \in \DU A m$,
 \[
	\eta^m_u(\alpha, n) = u(i)(n),
 \]
 where $i < i^*$ is the index with $\alpha = \alpha_i$. We set
 \[
	f_u = \BigCompCond [m < \omega] {T^m (\eta^m_u)}.
 \]
 Now, if $f_u$ was as required in \Lemma {Bound} and $T^m (\eta^m_u)
\Wkr T^{m +1} (\eta^{m +1}_u)$ for every $m < \omega$, then it would
follow, by the same lemma, that there is some $u$ and $r \in P^*$ such
that $f_u \Subset r$ and $r$ forces \Par {\Name \LadSys (\delta^*) =
\Ladder y} and $(\NameCol [_{\delta^*, n}] b \not= 0)$ for all $n <
\omega$. By the definition of $f_u$, $r$ would be a lower bound of the
descending chain \Seq {T^m (\eta^m_u)} {m < \omega} of conditions. So
to prove the claim it suffices to check that the conditions $T^m
(\eta^m_u)$, $m < \omega$, form a descending chain of conditions and
$f_u$ satisfies the properties wanted in \Lemma {Bound}.

\LemmaItem {Bound} {union}
 The function $f_u$ is well-defined since for all $i$ and $n$ such
that $(\alpha_i, n) \in \DU A m$,
 \[
	\eta^m_u (\alpha_i, n) = u(i)(n)
	= \eta^{m+1}_u (\alpha_i, n),
 \]
 i.e., $\eta^m_u = \Res {\eta^{m +1}_u} {\DU A m}$, and so by
\ItemOfDefinition {TreeSystem} {stronger}, $T^m (\eta^m_u) \geq T^{m
+1} (\eta^{m +1}_u)$. For each \Function u {i^*} {\Functions \omega
F}, $\Dom {f_u} \Subset \Braces 0 \Union A^*$ by \ItemOfDefinition
{TreeSystem} {domain}, and $f_u$ has height $\delta^*$ by
\ItemOfDefinition {TreeSystem} {height}.

\LemmaItem {Bound} {tree}
 Suppose \Function {u, v} {i^*} {\Functions \omega F}, $0 < i < i^*$,
and \EqualRes u v i. For all $m < \omega$ and $(\alpha, n) \in \DU A m
\Inter \alpha_i$, $\alpha$ must be $\alpha_j$ for some $j < i$ since
$\alpha < \alpha_i$, and furthermore,
 \[
	\eta^m_u (\alpha_j, n)
	= u(j)(n) = v(j)(n)
	= \eta^m_v (\alpha_j, n).
 \]
 Thus for each $m < \omega$, \EqualRes {\eta^m_u} {\eta^m_v}
{\alpha_i}, and by \ItemOfDefinition {TreeSystem} {Tree}, \EqualRes
{T^m (\eta^m_u)} {T^m (\eta^m_v)} {\alpha_i}. Consequently, for all
$\beta \in \Dom {f_u} \Inter \alpha_i = \Dom {f_v} \Inter \alpha_i$,
 \[
	\Forces [\beta] f_u(\beta)
	= \BigUnion [m < \omega] T^m (\eta^m_u) (\beta)
	= \BigUnion [m < \omega] T^m (\eta^m_v) (\beta)
	= f_v(\beta),
 \]
 and we may assume $f_u(\beta)$ is the same name as $f_v(\beta)$,
i.e., \EqualRes {f_u} {f_v} {\alpha_i}.

\LemmaItem {Bound} {almost}
 Let \Function u {i^*} {\Functions \omega F} and $i < i^*$ be such
that $\alpha_i \in \Dom {f_u}$. Then $\alpha_i$ is not $T^m
(\eta^m_u)$-trivial for any $m < \omega$. Let $n < \omega$ be such
that $\alpha_i \in A_n$. Then for each $m \geq n$, $\alpha_i \in A_m$,
and by \ItemOfDefinition {TreeSystem} {m_value},
 \ARRAY[lll]{
	    f_u (\alpha_i) (y_m)
	&=& T^m (\eta^m_u) (\alpha_i) (y_m) \\
	&=& \eta^m_u (\alpha_i, m) \\
	&=&  u(i)(m).
 }
\end{Proof}%
%

\end{LEMMA}


Now the main problem to be solved is the existence of a tree system
where each condition tree decides enough information about the
uniformizing function $\Name h$.

\begin{LEMMA} {TreeSystem}
 There exist a countable subset $A^*$ of $\alpha^* \Minus \Braces 0$,
a filtration \Filtration A of $A^*$, $\delta^* \in S$, an increasing
sequence $\bar \epsilon$ of ordinals with limit $\delta^*$, a ladder
$\Ladder y$ on $\delta^*$, and an $(\bar \epsilon, \Ladder y)$-tree
system \TreeSystem T on \Filtration A such that for all $m < \omega$
and $\eta \in \Ind A m$, $T^m_\eta \Forces [\alpha^*] \Name h(y_m) =
0$.
 \end{LEMMA}

We get the desired contradiction using the tree system given by the
preceding lemma together with \Lemma {TreeSystemBound}. Namely, a
lower bound $r \in P^*$ given by \Lemma {TreeSystemBound} satisfies
 \[
	r \Forces [\alpha^*]
	\Name [_{\delta^*,m}] \LadSys = y_m \And
	\NameCol [_{\delta,m}] b \not= 0, \ForAll m < \omega.
 \]
 On the other hand, \Lemma {TreeSystem} ensures that the lower bound
$r$ also satisfies the following condition:
 \[
	r \Forces [\alpha^*] \Name h(y_m) = 0, \ForAll m < \omega.
 \]
 It follows that $r \Str_{\alpha^*} p^*$, $\delta^* \in S$, and $r
\Forces [\alpha^*] \Par {\Name h(\Name \LadSys (\delta^*))
\not\IsAlmost \NameCol b (\delta^*)}$ contrary to our assumption
\PropertyOnPage {contra}. So, to prove \Theorem {Colouring} it
suffices to show that \Lemma {TreeSystem} holds. To achieve this goal
we have to analyze the relation between the values of conditions and
the value of \Name h in detail. Therefore we shall delay the proof of
\Lemma {TreeSystem} until the end of this subsection.



The following is a strengthening of \Lemma {not_decide}.

\begin{LEMMA} {not_decide(H)}
 Suppose $\alpha < \alpha^*$, $d \in F$, $Y$ is an unbounded subset of
\Vec, $p \in P^*$, and $H$ is a $P_\alpha$-generic set over $V$
containing \Res p \alpha. Then there is an unbounded subset $Z$ of $Y$
and for every $z \in Z$ a condition $q^z \in P^*$ satisfying
 \ARRAY{
	q^z \Str_{\alpha^*} p, \\
	\Res {q^z} \alpha \in H, \\
	q^z \Forces [\alpha^*] \Name h(z) = d.
 }

\begin{Proof}
 Suppose the lemma fails, and fix $\alpha, p, d, Y$, and $H$. Recall
what \Fact {Quotient} asserts and note that in $V[H]$ the condition
$p$ belongs to $P_{\alpha, \alpha^*}$. Consider the set $Y$ and $p$ in
$V[H]$. By our assumption, for all unbounded $Z \Subset Y$ there must
be some $z \in Z$ such that
 \[
	\ForAll s \in P^*,
	\If \Res s \alpha \in H \And s \Str_{\alpha^*} p
	\Then s \DoesNotForce [\alpha^*] \Name h(z) = d.
 \]
 Directly by \ItemOfFact {Quotient} {ab->b}, the following holds in
$V[H]$,
\[
	\ForAll r \in P_{\alpha,\alpha^*},
	\If r \Str_{\alpha,\alpha^*} p \Then
	r \DoesNotForce [\alpha,\alpha^*] \NName h(z) = d.
\]
Hence, for all sets $Z_\theta = \Set [\big] {y \in Y} {\theta <
\min(\Supp y)}$, where $\theta < \omega_1$, there is $z_\theta \in
Z_\theta$ such that in $V[H]$, for every $r \Str_{\alpha, \alpha^*} p$
in $P_{\alpha, \alpha^*}$ there is a condition $t \Str_{\alpha,
\alpha^*} r$ in $P_{\alpha, \alpha^*}$ for which $t \Forces [\alpha,
\alpha^*] \NName h(z_\theta) \not= d$. This means that in $V[H]$ the
collection of those conditions which forces $(\NName h(z_\theta) \not=
d)$ is dense below $p$ in the sense of $P_{\alpha, \alpha^*}$. Thus in
$V[H]$, $p \Forces [\alpha, \alpha^*] (\NName h(z_\theta) \not= d)$
for all $\theta < \omega_1$. By \ItemOfFact {Quotient} {ab->b} there
is $s \Str_{\alpha^*} p$ in $P^*$ forcing $(\Name h(z_\theta) \not=
d)$, for all $\theta < \omega_1$. This contradicts \Lemma
{not_decide}.
\end{Proof}
\end{LEMMA}

\begin{DEFINITION} {Pos}
 For all nonzero $\alpha < \alpha^*$ and $p \in P^*$ we define \Pos
\alpha p to be the set of tuples $(c_0, d_0, c_1, d_1) \in F^4$
satisfying the following requirement. There is an unbounded subset $Y$
of \Vec, and for each $y \in Y$ conditions $q_i^y \Str_{\alpha^*} p$
in $P^*$, $i = 0, 1$, such that
 \begin{ITEMS}

\ITEM {EqualRes} \EqualRes {q_0^y} {q_1^y} \alpha;

\ITEM {cond} either $\alpha$ is both $q_0^y$-trivial and
$q_1^y$-trivial, or otherwise $q_i^y(\alpha)(y) = c_i$ for both $i = 0
\And 1$;

\ITEM {h} $q_i^y \Forces [\alpha^*] \Name h(y) = d_i$ for both $i = 0
\And 1$.

 \end{ITEMS}
 \end{DEFINITION}

In the following lemma, the property \LemmaItem {Pos} {main} will be
the principal one later on.

\begin{LEMMA} {Pos}
 \begin{ITEMS}

\ITEM {trivial}
 If $p \in P^*$ and nonzero $\alpha < \alpha^*$ are such that there is
$q \Str_{\alpha^*} p$ in $P^*$ for which $\alpha$ is $q$-trivial, then
$(c, d_0, c, d_1) \in \Pos \alpha p$ for all $c, d_0, d_1 \in F$.

\ITEM {implication}
 If $\alpha < \alpha^*$ nonzero, $p \in P^*$, and $(c_0, d, c_1, d)
\in \Pos \alpha p$, where $c_0 \not= c_1, d \in F$, then there are $c,
d_0 \not= d_1 \in F$ such that $(c, d_0, c, d_1) \in \Pos \alpha p$.

\ITEM {main}
 For all $p \in P^*$ and nonzero $\alpha < \alpha^*$, there are $c,
d_0 \not= d_1 \in F$ such that $(c, d_0, c, d_1) \in \Pos \alpha p$.

\end{ITEMS}



\begin{Proof}%
\ProofOfItem {trivial}
 Let $H$ be a $P_\alpha$-generic set over $V$ containing \Res q
\alpha. By \Lemma {not_decide(H)} there are an unbounded subset $Y$ of
\Vec and conditions \Seq {q_0^y} {y \in Y} in $P^*$ such that for
every $y \in Y$,
 \ARRAY{
	q_0^y \Str q, \\
	\Res {q_0^y} \alpha \in H, \\
	q_0^y \Forces [\alpha^*] \Name h(y) = d_0.
 }
 By the same lemma there are an unbounded subset $Z$ of $Y$ and
conditions \Seq {q_1^y} {y \in Z} in $P^*$ such that
 \ARRAY{
	q_1^y \Str q, \\
	\Res {q_1^y} \alpha \in H, \\
	q_1^y \Forces [\alpha^*] \Name h(y) = d_1.
 }
 By \Fact {Res} there are, for $y \in Z$ and $i = 0, 1$, $r_i^y \Str
q_i^y$ in $P^*$ such that \EqualRes {r_0^y} {r_1^y} \alpha. Then for
all $c \in F$, the unbounded subset $Z$ of \Vec and the conditions
\Seq {r_i^y} {i = 0, 1 \And y \in Z} exemplify that $(c, d_0, c, d_1)
\in \Pos \alpha q \Subset \Pos \alpha p$. Observe that $\alpha$ is
$r_i^y$-trivial for $i = 0, 1$.

For the rest of the proof, we can restrict ourselves to the case that
\Res p \alpha forces \Name [_\alpha] Q to be nontrivial by \Item
{trivial}.

\ProofOfItem {implication}
 Suppose an unbounded subset $Y$ of \Vec and conditions $q_0^y, q_1^y
\Str p$ for $y \in Y$ exemplify that $(c_0, d, c_1, d) \in \Pos \alpha
p$. By the nontriviality of $\alpha$ we assume that for $i = 0, 1$ and
$y \in Y$,
 \ARRAY{ 
	\EqualRes {q_0^y} {q_1^y} \alpha, \\
	q_i^y(\alpha)(y) = c_i, \\
	q_i^y \Forces [\alpha^*] \Name h(y) = d.
 }

Consider some $y \in Y$ and $q_0^y$. Let $H$ be a $P_\alpha$-generic
set over $V$ such that $\EqualRes {q_0^y} {q_1^y} \alpha \in H$. By
\Lemma {not_decide(H)} there must be an unbounded subset $Z_0^y$ of
\Vec satisfying for all $z \in Z_0^y$ that $\max (\Supp y) < \min (
\Supp z)$ and there is $r_0^{y, z} \in P^*$ such that $r_0^{y, z} \Str
q_0^y$, $\Res {r_0^{y, z}} \alpha \in H$, and $r_0^{y, z} \Forces
[\alpha^*] \Name h(z) = 0$. Since $Z_0^y$ is unbounded, we can use the
same lemma again. Hence there must be some $z^y \in Z_0^y$ and a
condition $r_1^{y, z^y} \Str q_1^y$ in $P^*$ such that $\Res {r_1^{y,
z^y}} \alpha \in H$, and $r_1^{y, z^y} \Forces [\alpha^*] \Name h(z^y)
= 1$. By \Fact {Res} there are in $P^*$ conditions $s_i^y \Str r_i^{y,
z^y}$ for $i = 0, 1$ such that \EqualRes {s_0^y} {s_1^y} \alpha.

By \ItemOfLemma {Extension} {P_alpha}, we may assume that $\Dom {z^y}
\Subset \Dom {s_i^y (\alpha)}$ for both $i = 0 \And 1$. Since $F$ is
countable and $Y$ is uncountable, there is an unbounded subset $Z$ of
$Y$ and $(a_0, a_1) \in F^2$ such that the pair $(s_0^y (\alpha)
(z^y), s_1^y (\alpha) (z^y))$ is $(a_0, a_1)$ for every $y \in Z$.

Define $e_0 = a_1 - a_0$ and $e_1 = c_0 - c_1$. Since $c_0 \not= c_1$,
$e_1$ is not 0 \Note {$e_0$ might be 0}. Now, for all $i = 0, 1$ and
$y \in Z$ the following hold
 \ARRAY{
	s_i^y(\alpha)(e_0 y + e_1 z^y) = e_0 c_i + e_1 a_i, \\
	s_i^y \Forces [\alpha^*] \Name h(e_0 y + e_1 z^y) = e_0 d + e_1 i.
 }
 Consequently, the unbounded subset \Set {(e_0 y + e_1 z^y)} {y \in Z}
of \Vec and the conditions $s_i^y$, for $i=0, 1$ and $y \in Z$,
exemplify that $(c, d_0, c, d_1) \in \Pos \alpha p$, where $c = e_0
c_0 + e_1 a_0$ \Note {$= e_0 c_1 + e_1 a_1 $}, $d_0 = e_0 d + e_1 0$,
and $d_1 = e_0 d + e_1 1$. Clearly, $d_0 \not= d_1$.

\ProofOfItem {main}
 We may assume that \Res p \alpha decides the value of \Dom {p
(\alpha)}. Suppose, contrary to the claim, that there are no elements
$c, d_0 \not= d_1$ in $F$ such that $(c, d_0, c, d_1) \in \Pos \alpha
p$. By \Item {implication} this implies that there are no $c_0 \not=
c_1, d \in F$ satisfying $(c_0, d, c_1, d) \in \Pos \alpha p$ either.

Let $H$ be a $P_\alpha$-generic set over $V$ such that $\Res p \alpha
\in H$. Define $\SPos H p$ to be the set of all $(\xi, c, d) \in
\omega_1 \times F \times F$ such that there is $q \in P^*$ satisfying
the following requirements:
 \ARRAY{
	q \Str_{\alpha^*} p,	\\
	\Res q \alpha \in H,	\\
	q(\alpha)(\xi) = c,	\\
	q \Forces [\alpha^*] \Name h(\xi) = d.
 }
 It is easy to see, using \Fact {Quotient}, that for all $\xi <
\omega_1$ satisfying $\xi \not\in \Dom {p(\alpha)}$, and $c \in F$,
there is $d \in F$ such that $(\xi, c, d) \in \SPos H p$. Namely, by
\ItemOfLemma {Extension} {P_alpha} there is $q \Str p$ for which $q
(\alpha) (\xi) = c$ and $\EqualRes q p \alpha \in H$. Since $\Res q
\alpha \in H$, and $q \Forces [\alpha^*] (\Function {\Name h} {\SName
[_1] \omega} F)$, the following holds in $V[H]$ by \ItemOfFact
{Quotient} {b->ab}: there are $r \Str q$ in $P_{\alpha, \alpha^*}$ and
$d \in F$ for which $r \Forces [\alpha, \alpha^*] \NName h(\xi) =
d$. By \ItemOfFact {Quotient} {ab->b} there is $s \Str r$ in $P^*$
satisfying $\Res s \alpha \in H$ and $s \Forces [\alpha^*] \Name
h(\xi) = d$. So, $s$ exemplifies $(\xi, c, d) \in \SPos H p$.

Another easy property is that if there is an unbounded subset $I$ of
$\omega_1$ and $c_0, c_1, d_0, d_1 \in F$ such that for every $\xi \in
I$ both $(\xi, c_0, d_0)$ and $(\xi, c_1, d_1)$ are in $\SPos H p$,
then $(c_0, d_0, c_1, d_1)$ is in \Pos \alpha p. Namely, if for $\xi
\in I$ the conditions $q_i^\xi \Str p$, $i = 0,1$, exemplify that
$(\xi, c_i, d_i) \in \SPos H p$, then both \Res {q_0^\xi} \alpha and
\Res {q_1^\xi} \alpha belong to $H$. By \Fact {Res} there are $r_i^\xi
\Str q_i^\xi$ in $P^*$, for $i = 0, 1$ and $\xi \in I$, such that
$\EqualRes {r_0^\xi} {r_1^\xi} \alpha$. The set \Set {\Gen \xi} {\xi
\in I} and the conditions $r_i^\xi$, for $i = 0, 1$ and $\xi \in I$,
exemplify that $(c_0, d_0, c_1, d_1) \in \Pos \alpha p$. Observe that
these two simple observations together imply that \Pos \alpha p is
always nonempty.

It follows from our initial assumptions that we can fix $\mu' <
\omega_1$ such that the definition
 \[
	\pi_\xi(c) = d \Iff (\xi, c, d) \in \SPos H p
 \]
 yields in $V[H]$ an injective function \Function {\pi_\xi} F F when
$\mu' \leq \xi < \omega_1$. Since $F$ is finite each $\pi_\xi$ is in
fact a permutation of $F$. From the definition of $\SPos H p$ it
follows that $p \Forces [\alpha^*] \Par {\pi_\xi(\Name [_\alpha]
g(\xi)) = \Name h (\xi)}$ for all $\mu' \leq \xi < \omega_1$.

A function \Function \psi F F is a \Def {line} if there are $k, m \in
F$ such that $\psi(a) = k a + m$ for all $a \in F$ \Note {$k$ is the
\Def {slope} of the line}.

Our proof of \Item {main} will have the following structure.
 \begin{enumerate}

\ITEM {not_line}
 First we assume that there are unboundedly many $\xi < \omega_1$ such
that $\pi_\xi$ is not a line. It will follow that there are $c, d_0
\not= d_1 \in F$ such that $(c, d_0, c, d_1) \in \Pos \alpha p$,
contrary to our initial assumption.

\ITEM {line}
 We assume the converse of \Item {not_line}, i.e., we suppose $\mu <
\omega_1$ is a limit such that $\mu' \leq \mu$ and for every $\mu \leq
\xi < \omega_1$,
 \TEXTPROPERTY{line}{
 $k_\xi$ and $m_\xi$ are elements in $F$ such that $\pi_\xi(a) = k_\xi
a + m_\xi$ holds for all $a \in F$ in $V[H]$.
 }%
 Since each $\pi_\xi$ is injective $k_\xi \not= 0$ for every $\mu \leq
\xi < \omega_1$. Using this assumption we shall make two more steps.
 \begin{enumerate}

 \ITEM {constant}
 We show that
 \TEXTPROPERTY{nonconstant}{
 there is no $\theta < \omega_1$ and $e \in F$ such that $k_\xi = e$
whenever $\Max {\theta, \mu} \leq \xi < \omega_1$.
 }%
 Observe that this is the only part of the proof of the theorem where
the condition $(\NameCol [^\alpha] a \not\in \GenBy {\NameCol b} +
\NameUnifSet)$ in \Definition {Q_alpha} is essential, i.e., that we do
not ``kill'' colourings which are too ``close'' to the generic
colouring \NameCol b.

 \ITEM {nonconstant}
 The last case is that for all $\xi \geq \mu$ there is $\zeta > \xi$
such that $k_\xi \not= k_\zeta$, i.e., the slopes of lines $\pi_\xi$,
$\pi_\zeta$ are different. This will yield that there are $c_0 \not=
c_1, d \in F$ such that $(c_0, d, c_1, d) \in \Pos \alpha p$, contrary
to our initial assumption.

 \end{enumerate}
\end{enumerate}

\ProofOfItem {not_line}
 We shall show that for each $\theta < \omega_1$ there are $y^\theta
\in \Vec$, conditions $q^\theta, r^\theta \Str p$ in $P^*$, and
elements $c^\theta, d^\theta \not= e^\theta$ in $F$ such that $\min
(\Supp {y^\theta}) > \theta$ and
 \ARRAY{
	\Res {q^\theta} \alpha = \Res {r^\theta} \alpha, \\
	q^\theta(\alpha)(y^\theta) = c^\theta
	= r^\theta(\alpha)(y^\theta), \\
	q^\theta \Forces [\alpha^*]
		\Name h(y^\theta) = d^\theta, \\
	r^\theta \Forces [\alpha^*]
		\Name h(y^\theta) = e^\theta. \\
 }
 Since the choice of $\theta$ will be arbitrary, it will follow that
there are uncountable $I \Subset \omega_1$ and $c, d \not= e \in F$
such that for every $\theta \in I$, $c^\theta = c$, $d^\theta = d$,
and $e^\theta = e$. Then the unbounded subset \Set {y^\theta} {\theta
\in I} of \Vec and conditions \Seq {q^\theta, r^\theta} {\theta \in I}
will exemplify that $(c, d, c, e)$ is in $\Pos \alpha p$, where $d
\not= e$, contrary to our initial assumption.

Let $\theta < \omega_1$ be given. Since there are uncountably many
$\xi < \omega_1$ for which $\pi_\xi$ is not a line and only finitely
many permutations of $F$, fix $\xi < \zeta < \omega_1$ such that $\Max
{\mu', \theta, \Dom {p(\alpha)}} < \xi$ and $\pi_\xi = \pi_\zeta$ is
not a line. Let $\pi$ be the function $\pi_\xi = \pi_\zeta$. Fix
arbitrary $a \not= b_0 \in F$, and let $\psi_0$ be the line satisfying
$\psi_0 (a) = \pi (a)$ and $\psi_0 (b_0) = \pi (b_0)$. Since $\pi$ is
not a line there is $b_1 \in F$ for which $\pi (b_1) \not= \psi_0
(b_1)$. Let $\psi_1$ be the line for which $\psi_1 (a) = \pi(a)$ and
$\psi_1 (b_1) = \pi (b_1)$.

By \ItemOfLemma {Extension} {P_alpha} and since $\Res p \alpha$ forces
\Par {\Name [_\alpha] Q = \NameUnifFunc [^\alpha] a}, there is a
condition $q^\theta \in P^*$ such that
 \ARRAY{
	\EqualRes p {q^\theta} \alpha, \\
	q^\theta \Str p, \\
	q^\theta(\alpha)(\xi) = a = q^\theta(\alpha)(\zeta).
 }
 By the same lemma again, there is $r^\theta \in P^*$ such that
 \ARRAY{
	\EqualRes p {r^\theta} \alpha, \\
	r^\theta \Str p, \\
	r^\theta(\alpha)(\xi) = b_0 \And r^\theta(\alpha)(\zeta) = b_1.
 }
 Hence $\EqualRes {q^\theta} {r^\theta} \alpha \in H$. From the
definition of $\pi_\xi$ and $\pi_\zeta$ it follows that
 \ARRAY{
	 q^\theta \Forces [\alpha^*]
		\Name h(\xi)
		= \pi_\xi(q^\theta(\alpha)(\xi))
		= \psi_0(a)
		\And
		\Name h(\zeta)
		= \pi_\zeta(q^\theta(\alpha)(\zeta))
		= \psi_1(a).
 }
 \Note {A proof of this fact is a reasoning concerning $\Forces
[\alpha^*]$ and $\Forces [\alpha, \alpha^*]$ similar to what we have
done many times earlier.} Analogously, $r^\theta$ satisfies
		$r^\theta \Forces [\alpha^*] \Par {
		     \Name h(\xi) = \psi_0(b_0)
		\And \Name h(\zeta) = \psi_1(b_1)}$.

Define $e_0 = b_1 - a$ and $e_1 = a - b_0$. Since $a \not= b_0$ and $a
\not= b_1$ both $e_0$ and $e_1$ are nonzero. Define $y^\theta = (e_0
\Gen \xi + e_1 \Gen \zeta)$ and $a^\theta = e_0 a + e_1 a$ \Note {$=
e_0 b_0 + e_1 b_1$}. Then
 \ARRAY{
	q^\theta(\alpha)(y^\theta) = e_0 a + e_1 a = a^\theta
	= e_0 b_0 + e_1 b_1 = r^\theta(\alpha)(y^\theta).
 }
 Moreover,
 \ARRAY{
	q^\theta \Forces [\alpha^*] \Name h(y^\theta)
	= e_0 \Name h(\xi) + e_1 \Name h(\zeta)
	= e_0 \psi_0(a) + e_1 \psi_1(a),
 }
 and
 \ARRAY{
	r^\theta \Forces [\alpha^*] \Name h(y^\theta)
	= e_0 \Name h(\xi) + e_1 \Name h(\zeta)
	= e_0 \psi_0(b_0) + e_1 \psi_1(b_1).
 }
 Define $d^\theta = e_0 \psi_0(a) + e_1 \psi_1(a)$ and $e^\theta = e_0
\psi_0(b_0) + e_1 \psi_1(b_1)$. Then $d^\theta \not=
e^\theta$. Namely, if they are equal then
 \[
	e_0 \psi_0(a) + e_1 \psi_1(a)
	= e_0 \psi_0(b_0) + e_1 \psi_1(b_1)
 \]
 implies
 \[
	e_0 k_0 (a - b_0) = e_1 k_1 (b_1 - a),
 \]
 where $k_0$ and $k_1$ are the slopes of the lines $\psi_0$ and
$\psi_1$ respectively \Note {i.e., for $i = 0, 1$ we assume $\psi_i
(a') = k_i a' + m_i$ for all $a' \in F$}. But from the choice of the
lines $\psi_i$ it follows that $k_0 \not= k_1$. Hence the preceding
equation contradicts our choice of $e_0$ and $e_1$.


\ProofOfItem {constant}
 Suppose $K$ is a $P_{\alpha, \alpha^*}$-generic set over $V[H]$
satisfying that $p \in K$ and for the elements $h = \Int K {\NName h}$
and $g_\alpha = \Int K {\NName [_\alpha] g}$, where the names \NName h
and \NName [_\alpha] g are given in \Fact {Quotient}, the equations
\Par {h(\xi) = \pi_\xi (g_\alpha (\xi))} for all $\mu \leq \xi <
\omega_1$ hold in $V[H][K]$.

 A proof of \Property {nonconstant} follows. Fix, contrary to the
claim, $\theta \geq \mu$ and $e$ satisfying \Property
{nonconstant}. Define in $V[H]$ a function \Function f {\omega_1} F
for all $\xi < \omega_1$ by
 \[
    \FunctionDefinition {f(\xi)}{
	\FunctionDefMidCase {0} {\If \xi < \theta}
	\FunctionDefOtherwise {\pi_\xi(0)}
    }
 \]
 Then $f$ satisfies in $V[H]$ the following equation for all $a \in F$
and $\theta \leq \xi < \omega_1$,
 \[
	f(\xi) = \pi_\xi(0) = m_\xi 
	= (e a + m_\xi) - e a = \pi_\xi(a) - e a.
 \]
 Hence, independently of what $g_\alpha$ is, the following equation
holds in $V[H][K]$ for all $\delta \in S$ and for almost all $n <
\omega$,
 \ARRAY[lll]{
	\Col [_{\delta,n}] b - e \Times \Col [_{\delta,n}^\alpha] a
	&=& h(\LadSys [_{\delta,n}]) 
		- e \Times g_\alpha(\LadSys [_{\delta,n}]) \\

	&=& \Par [\big] {\Sum [\xi < \delta]
		{e^{\delta, n}_\xi \Times h(\xi)} }
	  - e \Times \Par [\big] {\Sum [\xi < \delta]
		{e^{\delta, n}_\xi \Times g_\alpha(\xi)}
	} \\ 

	&=& \Sum [\xi < \delta] {
	    \Par [\Big] {
		e^{\delta,n}_\xi \Times
			\Par [\big] {h(\xi) - e \Times g_\alpha(\xi)}
	    }
	} \\ 
	&=& \Sum [\xi < \delta] {
	    \Par [\Big] {
		e^{\delta,n}_\xi \Times \Par [\big] {
			\pi_\xi(g_\alpha(\xi)) - e \Times g_\alpha(\xi)
		}
	    }
	} \\ 
	&=& \Sum [\xi < \delta] {e^{\delta,n}_\xi \Times f(\xi)},
}
 where each $\LadSys [_{\delta, n}]$ is assumed to be of the form \Sum
[\xi < \delta] {e^{\delta, n}_\xi \Gen \xi}.

But $f$ is already in $V[H]$. So, from \ItemOfLemma {Colouring}
{initial} it follows that $\Col b \Equivalent e \Times \Col [^\alpha]
a$, and hence, \Par {\Col [^\alpha] a \in \GenBy {\Col b} + \UnifSet}
holds in $V[H]$. By \Definition {Q_alpha}, \Int H {\Name [_\alpha] Q}
must be \Braces \One. Since $\Res p \alpha \in H$, this contradicts
our initial assumption that \Res p \alpha forces \Name [_\alpha] Q to
be nontrivial.


\ProofOfItem {nonconstant}
 If the size of $F$ is 2, then for every $\mu \leq \xi < \omega_1$ the
value of $k_\xi$ must be constantly 1 contradicting \Property
{nonconstant}. Hence the lemma holds if $F$ is of size $2$.

Now, $\Card F > 2$, \Property {line} holds, and $k_\xi \not= 0$ for
all $\mu < \xi < \omega_1$. Analogously to the case \Item {not_line},
to prove that there are $c \not= e, d \in F$ for which $(c, d, e, d)
\in \Pos \alpha p$, it suffices to show for arbitrary $\theta <
\omega_1$ the existence of $y^\theta \in \Vec$, and conditions
$q^\theta$, $r^\theta$ in $P^*$ satisfying
 \ARRAY[l]{

	\min (\Supp {y^\theta}) > \theta,  \\

	q^\theta, r^\theta \Str p, \\

	\Res {q^\theta} \alpha = \Res {r^\theta} \alpha, \\

	q^\theta(\alpha)(y^\theta) = c^\theta, \\

	r^\theta(\alpha)(y^\theta) = e^\theta, \\

	q^\theta \Forces [\alpha^*] \Name h(y^\theta) = d^\theta, \\

	r^\theta \Forces [\alpha^*] \Name h(y^\theta) = d^\theta.

 }

Let $\theta < \omega_1$ be given. Fix $\xi > \Max {\mu, \theta, \Dom
{p(\alpha)}}$ and $\zeta > \xi$ such that $k_\xi \not= k_\zeta$. As in
\Item {not_line} fix $q^\theta, r^\theta \Str p$ such that
 \ARRAY[l]{
	\EqualRes {q^\theta} {r^\theta} \alpha \in H, \\
	q^\theta(\alpha)(\xi) = 1 \And q^\theta(\alpha)(\zeta) = 1, \\
	r^\theta(\alpha)(\xi) = 2 \And r^\theta(\alpha)(\zeta) = 2.
 }
 Define $e_\xi = -k_\zeta$ and $e_\zeta = k_\xi$. Then $e_\xi k_\xi +
e_\zeta k_\zeta = 0$, and $e_\xi + e_\zeta \not= 0$ since $k_\xi \not=
k_\zeta$. If we let $y^\theta$ be $(e_\xi \Gen \xi + e_\zeta \Gen
\zeta)$, then
 \[
	q^\theta(\alpha)(y^\theta)
	= e_\xi \Times q^\theta(\alpha)(\xi)
	  + e_\zeta \Times q^\theta(\alpha)(\zeta)
	= e_\xi + e_\zeta,
 \]
 and
 \ARRAY[lll]{
	q^\theta \Forces [\alpha^*] \Name h(y^\theta)
	&=& e_\xi \Times (k_\xi + m_\xi)
	  + e_\zeta \Times (k_\zeta + m_\zeta) \\
	&=& (e_\xi k_\xi + e_\zeta k_\zeta)
	  + (e_\xi m_\xi + e_\zeta m_\zeta) \\
	&=& (e_\xi m_\xi + e_\zeta m_\zeta).
 }
 By a similar reasoning $r^\theta$ satisfies
 \ARRAY[lll]{
	r^\theta(\alpha)(y^\theta) &=& 2 (e_\xi + e_\zeta), \\
	r^\theta \Forces [\alpha^*] \Name h(y^\theta)
			&=& 2(e_\xi k_\xi + e_\zeta k_\zeta)
			  + (e_\xi m_\xi + e_\zeta m_\zeta)
			  = (e_\xi m_\xi + e_\zeta m_\zeta).
 }
 Hence $c^\theta = e_\xi + e_\zeta (\not= 0)$, $e^\theta = 2 c^\theta
(\not= c^\theta)$, and $d^\theta = e_\xi m_\xi + e_\zeta m_\zeta$ are
the desired elements of $F$.
\end{Proof}%
%

\end{LEMMA}



Now we can proceed with analyzing properties of condition
trees. Recall that \Filtration A is a filtration of $A^*$. Suppose $m
< \omega$, $T$ is an \DU A m-condition tree, $\eta \in \Ind A m$, and
$p$ is a condition in $P^*$ such that $\Res p \gamma \Str \Res
{T(\eta)} \gamma$ for $\gamma = \max A_m$. We define a function
\TreeComp T \eta p by setting for all $\nu \in \Ind A m$ that
\[
	\FunctionDefinition {\TreeComp T \eta p (\nu)}{
	\FunctionDefMidCase {p} {\If \nu = \eta}
	\FunctionDefMidCase {\CompCond {\Res p {\beta_\nu}} {T (\nu)}}
				{\Otherwise}
	}
\]
where $\beta_\nu = \Max [\EqualRes \nu \eta \gamma] {\gamma \in
A_m}$. Observe that for each $\nu \in \Ind A m$, $\TreeComp T \eta p
(\nu)$ is a condition in $P^*$ since $\Res p {\beta_\nu} \Str \Res
{T_\eta} {\beta_\nu} = \Res {T_\nu} {\beta_\nu}$. Hence, \TreeComp T
\eta p is an \DU A m-condition tree and $\TreeComp T \eta p \Str T$.

\begin{LEMMA} {Tree_height}
 Suppose $\epsilon < \omega_1$ and $T$ is an \DU A m-condition
tree. Then there is an \DU A m-condition tree $R \Str T$ of height
$\geq \epsilon$.

\begin{Proof}
 Suppose \Set {\eta_i} {i < k}, $k < \omega$, is an enumeration of
$\Ind A m$. We define by induction on $j \leq k$, \DU A m-condition
trees $R^j$ as follows. Let $R^0$ be $T$. Suppose $j < k$, $R^i$ for
all $i \leq j$ are defined, and the conditions $R^j (\eta_i)$, $i <
j$, are of height $\geq \epsilon$. By \ItemOfLemma {UnionAndHeight}
{height} there is $p \Str R^j (\eta_j)$ in $P^*$ having height greater
than $\epsilon$. We define $R^{j +1}$ to be \TreeComp {R^j} {\eta_j}
p. It follows that $R^k \Str T$ is an \DU A m-condition tree of height
$\geq \epsilon$.
 \end{Proof}
\end{LEMMA}

\begin{DEFINITION} {TPos}
 We fix the following notation for each $m < \omega$:
 \ARRAY [lll] {
	\Val A m &=& \Set {\tau} {\tau
	\Text{is a function from \Ind A m into} F}, \\

	\IInd A m &=& \Set {\Res \eta {\alpha+1}}
	{\alpha \in A_m \And \eta \in \Ind A m}, \\

	\IVal A m &=& \Set {\sigma} {\sigma
	\Text{is a function from \IInd A m into} F}.
 }
 Let $m < \omega$ and $T$ be an \DU A m-condition tree. For all $y \in
\Vec$ and $(\sigma, \tau) \in \IVal A m \times \Val A m$ we write
 \[
	T [y] \Follows (\sigma, \tau)
 \]
 if for each $\eta \in \Ind A m$, both of the requirements
 \[
	\ForEach \alpha \in A_m,
	\Text {either $\alpha$ is $T_\eta$-trivial or}
	T_\eta (\alpha) (y) = \sigma(\Res \eta {\alpha+1}),
 \]
 and
 \[
	T_\eta \Forces [\alpha^*] \Name h (y) = \tau(\eta),
 \]
 are satisfied. We define \TPos A m to be the set of all $(\sigma,
\tau) \in \IVal A m \times \Val A m$ with the following property. For
all \DU A m-condition trees $T$ there exist an unbounded subset $Y$ of
\Vec and for each $y \in Y$ an \DU A m-condition tree $T^y \Str T$
satisfying $T^y [y] \Follows (\sigma, \tau)$.
 \end{DEFINITION}

Suppose $m < \omega$ and $T$ is an \DU A m-condition tree. We set
 \ARRAY[lllll]{
	\Dec T &=& \big\{y \in \Vec &\SetSeparator&
	\ForAll \eta \in \Ind A m \And \alpha \in A_m, \\
	& & & & \alpha \Text{is $T_\eta$-trivial or}
	\Supp y \Subset \Dom {T_\eta (\alpha)} \big\}, \\
\\
	\NotDec T &=& \big\{y \in \Vec &\SetSeparator&
	\ForAll \eta \in \Ind A m \And \alpha \in A_m, \\
	& & & & \Supp y \not\Subset \Dom {T_\eta(\alpha)} \big\}, \\
\\
	\Dec [h] T &=& \big\{y \in \Vec &\SetSeparator&
	\ForEach \eta \in \Ind A m, \\
	& & & & T_\eta \Text{decides the value of} \Name h (y) \big\}.
 }
 For $i = 0, 1$, $(\sigma_i, \tau_i) \in \IVal A m \times \Val A m$
and $e_i \in F$ we define the sum
 \[
	e_0 \Times (\sigma_0, \tau_0) + e_1 \Times (\sigma_1, \tau_1)
 \]
 to be the pair $(\sigma, \tau) \in \IVal A m \times \Val A m$, where
for all $\upsilon \in \IInd A m$ and $\eta \in \Ind A m$
 \ARRAY[lll]{
	\sigma(\upsilon) &=&
	e_0 \Times \sigma_0 (\upsilon) + e_1 \Times \sigma_1 (\upsilon), \\
	\tau(\eta)   &=&
	e_0 \Times \tau_0 (\eta) + e_1 \Times \tau_1 (\eta).
 }

\begin{LEMMA} {Am-tree}
 Suppose $m < \omega$ and $T$ is an \DU A m-condition tree.
 \begin{ITEMS}

\ITEM {decides}
 For every $y \in \Vec$ there is an \DU A m-condition tree $R \Str T$
for which $y \in \Dec [inter] R$.

\ITEM {Follows}
 For all $y \in \NotDec T$ and $\sigma \in \IVal A m$ there are an \DU
A m-condition tree $R \Str T$ and $\tau \in \Val A m$ such that $R [y]
\Follows (\sigma, \tau)$.

\ITEM {Exists_tau}
 For every $\sigma \in \IVal A m$, there is $\tau \in \Val A m$ such
that $(\sigma, \tau) \in \TPos A m$.

\ITEM {addition}
 If $(\sigma_i, \tau_i) \in \TPos A m$ and $e_i \in F$, for $i = 0,
1$, then \Sum [i = 0, 1] {e_i \Times (\sigma_i, \tau_i)} is in \TPos A
m.

 \end{ITEMS}

\begin{Proof}
 \ProofOfItem {decides}
 Suppose $\Ind A m = \Set {\eta_i} {i < k}$. Let $R^0$ be $T$. Assume
\DU A m-condition trees $R^i$, $i \leq j < k$, are already defined.
 \TEXTPROPERTY{A}{
 By \ItemOfLemma {Colouring} {extension} there is $p \Str R^j
(\eta_j)$ in $P^*$ for which $\Supp y \Subset \Dom {p (\alpha)}$ for
all $\alpha \in A_m$.
 }%
 Assume $q \Str p$ in $P^*$ decides the value of $\Name h (y)$, and
define $R^{j +1}$ to be \TreeComp {R^j} {\eta_j} q. Then $y \in \Dec
[inter] {R^k}$.

\ProofOfItem {Follows}
 This is proved as \Item {decides}. The only difference is that
instead of \Property A the following is used:
 \begin{property}
 by \ItemOfLemma {Extension} {P_alpha} there is $p \Str R^j (\eta_j)$
in $P^*$ satisfying for each $\alpha \in A_m$ that either $\alpha$ is
$p$-trivial or otherwise $p (\alpha) (y) = \sigma (\Res {\eta_j}
{\alpha +1})$.
 \end{property}
 Then the function $\tau \in \Val A m$ satisfying $R^k [y] \Follows
(\sigma, \tau)$ is uniquely determined by $R^k$.

\ProofOfItem {Exists_tau}%
Since $T$ and the domains of the conditions in $T$ are countable there
must be a limit $\theta_T < \omega_1$ such that for every $y \in
\Vec$, $\min (\Supp y) > \theta_T$ implies $y \in \NotDec T$. Hence,
directly by \Item {Follows}, for every $y \in \NotDec T$ there are
$T^y \Str T$ and $\tau^y \in \Val A m$ satisfying $T^y [y] \Follows
(\sigma, \tau^y)$. Since \Val A m is countable and \NotDec T
uncountable, there must be an unbounded subset $Y$ of \NotDec T and
$\tau \in \Val A m$ such that $\tau = \tau^y$ for each $y \in Y$. Thus
$Y$ and the trees \Seq {T^y} {y \in Y} stronger than the arbitrary \DU
A m-condition tree $T$ exemplify $(\sigma, \tau) \in \TPos A m$.

\ProofOfItem {addition}
 Since $(\sigma_0, \tau_0) \in \TPos A m$ there are an unbounded
subset $Y$ of \Vec and for each $y \in Y$, an \DU A m-condition tree
$T_0^y \Str T$ satisfying $T_0^y[y] \Follows (\sigma_0,
\tau_0)$. Because $(\sigma_1, \tau_1) \in \TPos A m$, there exist for
each $y \in Y$ an \DU A m-condition tree $T_1^y \Str T_0^y$ and an
element $z_y \in \Vec$ such that $\max (\Supp y) < \min (\Supp {z_y})$
and $T_1^y [z_y] \Follows (\sigma_1, \tau_1)$. Consequently, for all
$y \in Y$,
 \[
	T^y_1 [e_0 y + e_1 z_y] \Follows
	e_0 \Times (\sigma_0, \tau_0) + e_1 \Times (\sigma_1, \tau_1).
 \]
 So the unbounded subset \Set {(e_0 y + e_1 z_y)} {y \in Y} of \Vec
and the trees \Seq {T_1^y} {y \in Y} stronger than an arbitrary $T$
exemplify that $\Sum [i = 0,1] {e_i \Times (\sigma_i, \tau_i)}$ is in
\TPos A m.
 \end{Proof}
\end{LEMMA} 

We let \IZero m be the $0$-function of \IVal A m and \Zero m be the
$0$-function of \Val A m. For all $\tau \in \Val A m$, $\eta \in \Ind
A m$, and $d \in F$, \ValMap \tau \eta d denotes the function in \Val
A m which is the same as $\tau$ except it maps $\eta$ into $d$.

\begin{LEMMA} {TPos}
 For every $\sigma' \in \IVal A m$ the pair $(\sigma', \Zero m)$ is in
\TPos A m.



\begin{Proof}%
We shall prove the following claim.
 \begin{property}
 For every $\eta_0 \in \Ind A m$ there are $(\sigma, \tau) \in \TPos A
m$ and $d_1 \in F$ such that $d_1 \not= \tau (\eta_0)$ and $(\sigma,
\ValMap \tau {\eta_0} {d_1})$ is in \TPos A m.
 \end {property}
 This suffices, because if the claim holds then by \ItemOfLemma
{Am-tree} {addition}
 \ARRAY[lll]{
    && \Div 1 {\tau(\eta_0) - d_1} \Times
    ((\sigma, \tau) - (\sigma, \ValMap \tau {\eta_0} {d_1})) \\
    &=& \Div 1 {\tau(\eta_0) - d_1} \Times
    (\IZero m, \ValMap {\Zero m} {\eta_0} {\tau(\eta_0) - d_1}) \\
    &=& (\IZero m, \ValMap {\Zero m} {\eta_0} 1) \in \TPos A m,
 }%
 for all $\eta_0 \in \Ind A m$. Furthermore, by \ItemOfLemma {Am-tree}
{Exists_tau}, there is $\tau' \in \Val A m$ for which $(\sigma',
\tau') \in \TPos A m$, and hence by \ItemOfLemma {Am-tree} {addition},
 \ARRAY[lll]{
    && (\sigma', \tau') -
    \Sum [\eta_0 \in \Ind A m]
    {\tau' (\eta_0) \Times (\IZero m, \ValMap {\Zero m} {\eta_0} 1)} \\
    &=& (\sigma', \tau') - (\IZero m, \tau') \\
    &=& (\sigma', \Zero m) \in \TPos A m.
 }

For the rest of the proof of the lemma let $\alpha$ be the maximal
element of $A_m$, $T$ be an \DU A m-condition tree, and $\eta_0$ be an
arbitrary element of \Ind A m. By \ItemOfLemma {Pos} {main} there are
$c, d_0 \not= d_1 \in F$, an unbounded subset $Z$ of \Vec, and
conditions $p_0^y, p_1^y \Str T (\eta_0)$, for each $y \in Z$,
exemplifying $(c, d_0, c, d_1) \in \Pos \alpha {T(\eta_0)}$. This
means that for all $y \in Z$, $i = 0, 1$, and $\beta \in A_m$,
 \ARRAY{
	\EqualRes {p_0^y} {p_1^y} \alpha, \\
	p_i^y \Forces [\alpha^*] \Name h (y) = d_i, \\
	p_0^y (\beta) (y) = p_1^y (\beta) (y) \Or
	\beta \Text{is $p_i^y$-trivial for both} i = 0 \And 1.
 }

By \ItemOfLemma {Am-tree} {decides} there is an \DU A m-condition tree
$T^y \Str \TreeComp T {\eta_0} {p_0^y}$ for every $y \in Z$ such that
$y \in \Dec [inter] {T^y}$. Since $Z$ is uncountable there must be an
unbounded subset $Y$ of $Z$ and $(\sigma, \tau) \in \IVal A m \times
\Val A m$ such that $T^y [y] \Follows (\sigma, \tau)$ for all $y \in
Y$. So $Y$ and the trees \Seq {T^y} {y \in Y} stronger than an
arbitrary tree $T$ exemplify $(\sigma, \tau)$ is in \TPos A m. Observe
that $T^y (\eta_0) \Str p_0^y$ implies $T^y (\eta_0) \Forces
[\alpha^*] \Name h (y) = \tau (\eta_0) = d_0$.

Now, the function
 \[
	R^y =
	\TreeComp {T^y} {\eta_0}
	{\CompCond {(\Res {T^y(\eta_0)} \alpha)} {p_1^y}}
 \]
 is a \DU A m-condition tree for each $y \in Y$, since $\Res {T^y
(\eta_0)} \alpha \Str \EqualRes {p_0^y} {p_1^y} \alpha$. Hence $Y$ and
\Seq {R^y} {y \in Y} exemplify $(\sigma, \ValMap \tau {\eta_0} {d_1})$
is in \TPos A m.
\end{Proof}%
%

\end{LEMMA}


We are now ready to give the last missing piece.

\begin{RefArea}{Env}{Lemma}{TreeSystem}
\begin{SeparateProof}
 Fix a countable elementary submodel \EM of $\HModel {\alpha^*} {}
(p^*, \Name h)$ satisfying $\EM \Inter {\omega_1} = \delta^* \in
S$. We define $A^* = \EM \Inter \alpha^*$. Let \Filtration A be a
filtration of $A^*$. Since the sets $A_m \Subset A^* \Subset \EM$, $m
< \omega$, are finite they belong to \EM as well as the sets \Ind A m,
\IInd A m, and \Val A m. Let $\bar \epsilon = \Seq {\epsilon_m} {m <
\omega}$ be an increasing sequence of ordinals with limit $\delta^*$.

For each $m < \omega$ we define the $A_m$-complete element of \IVal A
m to be the unique $\sigma \in \IVal A m$ for which $\sigma (\Res \eta
{\alpha +1}) = \eta (\alpha, m)$ for all $\eta \in \Ind A m$ and
$\alpha \in A_m$.

We define a ladder $\Ladder y = \Seq {y_m} {m < \omega}$ on $\delta^*$
and an $(\bar \epsilon, \Ladder y)$-tree system \Seq {T^m} {m <
\omega} on \Filtration A by induction on $m < \omega$. Our main tool
is \Lemma {TPos} which will ensure that $T^m_\eta$ forces \Par {\Name
h (y_m) = 0} for all $m < \omega$ and $\eta \in \Ind A m$. During the
induction we work inside \EM.

Suppose $m = 0$. We define a trivial \DU A 0-condition tree $R$ in \EM
by, $R (\eta) = p^*$ for each $\eta \in \Ind A 0$. Note that $\Dom
{p^*} \Subset \Braces 0 \Union A^*$. By \Lemma {Tree_height} there is
in \EM an \DU A 0-condition tree $R' \Str R$ which is of height $\geq
\epsilon_0$. By \Lemma {TPos} there are $y_0 \in \Vec \Inter \EM$ and
an \DU A 0-condition tree $T^0 \Str R'$ in \EM satisfying
 \[
	\epsilon_0 < \min(\Supp {y_0}) \And
	T^0 [y_0] \Follows (\sigma, \Zero m),
 \]
 where $\sigma$ is the $A_0$-complete element of \IVal A 0.

Similarly, when $y_m \in \Vec \Inter \EM$ and $T^m$ in \EM are already
defined, we can find $y_{m +1} \in \Vec \Inter \EM$ and an \DU A {m
+1}-condition tree $T^{m +1} \Str T^m$ in \EM satisfying
 \ARRAY{
	\Max {\epsilon_{m+1}, \max(\Supp {y_m})} < 
	\min (\Supp {y_{m+1}}), \\
	T^{m+1} \Text{is of height $\geq \epsilon_{m+1}$}, \\
	T^{m+1} [y_{m+1}] \Follows (\sigma, \Zero {m+1}),
 }
 where $\sigma \in \IVal A {m +1}$ is $A_{m +1}$-complete.

It follows directly from the definition above that \Ladder y is a
ladder on $\delta^*$ and for every $m < \omega$,
 \ARRAY{
	T^m \Text{is an \DU A m-condition tree}, \\
	\ForAll \eta \in \Ind A m,
		\Dom{T^m_\eta} \Subset \Braces 0 \Union A^*, \\
	T^m \Text{is of height $\geq \epsilon_m$ and $< \delta^*$}, \\
	T^{m+1} \Str T^m.
 }
Moreover, for each $m < \omega$ and $\eta \in \Ind A m$ the property
$T^m [y_m] \Follows (\sigma, \Zero m)$ guarantees that
 \[
	T^m_\eta \Forces [\alpha^*] \Name h (y_m) = \Zero m (\eta) = 0,
 \]
 and since $\sigma$ is $A_m$-complete,
 \begin{property}
	$\alpha$ is $T^m_\eta$-trivial or
	$T^m_\eta (\alpha) (y_m)
	= \sigma(\Res \eta {\alpha+1})
	= \eta (\alpha, m)$, for all $\alpha \in A_m$.
 \end{property}
 \end{SeparateProof}
\end{RefArea} 


\end{SUBSECTION}

\begin{SUBSECTION} {-} {Remarks} {Remarks}




There is a forcing notion which gives the conclusion of \Theorem
{Colouring} for all finite fields simultaneously. Namely, we defined
an iterated forcing $P_k = \CountLim \Seq {P_\alpha, \Name [^k_\alpha]
Q} {\alpha < \omega_2}$ for fixed $k$. The extended result would
follow if each \Name [^k_\alpha] Q was replaced by $\Name [^2_\alpha]
Q \times \Name [^3_\alpha] Q \times \dots$ where \Name [^i_\alpha] Q
takes care of the case $\pi(i) = (p, m)$ and $\pi$ is a coding for the
pairs of primes and positive integers. So $F_i$ would be the field of
size $p^m$ where $\pi(i) = (p, m)$. For example, to prove that for
each ``coordinate'' $i$ the cardinality of \Quotient {\ColSet [S,
F_i]} {\UnifSet [\LadSys, D]} is as wanted, it would suffice to
concentrate on one coordinate $i$, and define the condition trees and
systems, \Pos \alpha p, etc., only for fixed $i$. Hence an assumption
that the size is wrong for some $i$ would lead to a contradiction in
the same way as in \Subsection {Nonuniform}.

It is possible to have a \Vec [F]-ladder system on $S$ such that
$\Card {\Quotient {\ColSet [S, F]} {\UnifSet [\LadSys, D]}} =
\aleph_0$. A proof of this fact would be a forcing argument just like
the one we have given. The only difference is that instead of one
generic colouring \NameCol b, one should add generic colourings \Seq
{\NameCol [_m] b} {m < \omega} by defining $Q_0 = \InitLadSet \times
\Functions \omega \InitColSet$. Then by replacing $\GenBy [F] {\Col b}
+ \UnifSet$ with
	$(\GenBy [F] {\Col [_0] b, \Col [_1] b, \dots} + \UnifSet)$
the desired result would follow. The conclusion of such a generalized
theorem would be
 $\Forces [P]
	\Card [\big] {
		\Quotient {\NameColSet [S, F]}
			  {\NameUnifSet [\Name \LadSys, D]}
	} =
	\Card [\big] {
		\GenBy [F] {\NameCol [_0] b, \NameCol [_1] b, \dots}
	} = \aleph_0$.
 Other changes would be, for example, that \Lemma {Nonuniform} would
have the form $\Forces [1] \Quote {\If \chi \in \GenBy [F] {\NameCol
[_0] b, \NameCol [_1] b, \dots} \Then \chi \not\in \UnifSet}$, and
analogous changes would be needed in \Lemma {Pos}.

\begin{ShComment}%
We may also continue the iteration longer than $\omega_2$ and get the
consistency of our main result with \CH + ``any reasonable value for
\AlephExp 1''. The $\aleph_2$-c.c. for such a forcing follows from the
use of pic \cite {ProperForcing} or better \cite [Section 2 of Chapter
8] {Sh:f}.

During the given proof, for example in \Lemma {Distributive}, it is
possible to use the general claim on preservation of $(\omega_1 \Minus
S)$-complete forcing notions and the preservation of properness for
the preservation of stationarity \cite [Chapter 5] {ProperForcing} or
\cite [Chapter 5] {Sh:f}. But this does not, however, help with the
main problem.
\end{ShComment}


\end{SUBSECTION}


\end{SECTION}

\begin{SECTION} {+} {Models} {The Models}




As in the preceding sections, we assume that $S \Subset \omega_1$ is
a set of limit ordinals, $F$ is a field, $D$ is a filter over $\omega$
including all cofinite sets of $\omega$, \Vec is the vector space over
$F$ freely generated by \Seq {\Gen \xi} {\xi <\omega_1}, \LadSys is a
\Vec-ladder system on $S$, \ColSet denotes the set of all
$F$-colourings on $S$, and \UnifSet is the set of all uniform
colourings.

Let \M be a model of vocabulary \Voc, $0 < n < \omega$, and $R \in
\Voc$ a relation symbol with $n +1$ many places. We say that $R$ is a
partial function in \M if there are $X \Subset \M^n$ and $Y \Subset
\M$ such that the interpretation $R^\M$ of the symbol $R$ in \M is a
function from $X$ into $Y$. For all relations $R \in \Voc$, which are
partial functions in \M, $R^\M(x) = y$ means $\ConCat x {\SimpleSeq y}
\in R^\M$, and atomic formulas $R(x, y)$ are written in the form $R(x)
= y$.

\begin{DEFINITION}{models}
 We define a vocabulary \Voc and for all $\Col a \in \ColSet$ models
\MC a of vocabulary \Voc by the following stipulations:
 \begin{ITEMS}

\ITEM {domain}
 Each model \MC a has the same domain $(S \times \ODFun \omega F)
\Union (\Vec \times F)$, where
 \[
	\ODFun \omega F =
	\Set [\big] {u \in \Functions \omega F} {
		\Set {n \in \omega} {u(n) = 0} \in D
	}.
 \]

\ITEM {R_VxF}
 For each $y \in \Vec$, $R_y$ is a unary relation symbol in \Voc and
${R_y}^{\MC a} = \Braces y \times F$.

\ITEM {R_SxO}
 For each $\delta \in S$, $R_\delta$ is a unary relation symbol in
\Voc and ${R_\delta}^{\MC a} = \Braces \delta \times \ODFun \omega F$.

\ITEM {Pr}
 For each $n < \omega$, \Pr n a denotes a function from $S \times
\ODFun \omega F$ into $\Vec \times F$ defined for all $(\delta, u) \in
S \times \ODFun \omega F$ by
 \[
	\Pr n a (\delta, u)
	= \Par {\LadSys [_{\delta,n}], \Col [_{\delta,n}] a +_F u(n)}.
 \]
 For each $n < \omega$, $\Pr [-] n {}$ is a binary relation in \Voc
and ${\Pr [-] n {}}^{\MC a} = \Pr n a$. So \Pr [-] n {} is a partial
function in \MC a.

\ITEM {+b}
 For all $b \in F$, $\P b \in \Voc$, \Function {{\P b}^{\MC a}} {\Vec
\times F} {\Vec \times F}, and for all $(y, c) \in \Vec \times F$,
 \[
	{\P b}^{\MC a} (y, c) = (y, c +_F b).
 \]

\ITEM {+u}
 For all $u \in \ODFun \omega F$, $\P u \in \Voc$, \Function {{\P
u}^{\MC a}} {S \times \ODFun \omega F} {S \times \ODFun \omega F}, and
for all $(\delta, v) \in S \times \ODFun \omega F$,
 \[
	{\P u}^{\MC a} (\delta, v)
	= (\delta, v +_{(\ODFun \omega F)} u),
 \]
 where $v +_{(\ODFun \omega F)} u$ is the function in \ODFun \omega F
defined for all $n < \omega$ by $(v +_{(\ODFun \omega F)} u) (n) =
v(n) +_F u(n)$.

\ITEM {+}
 The symbol $+$ is in \Voc, \Function {+^{\MC a}} {{(\Vec \times
F)}^2} {\Vec \times F}, and for all $(y, b), (z, c) \in \Vec \times
F$,
 \[
	(y,b) +^{\MC a} (z,c) = (y +_{\Vec} z, b +_F c).
 \]

\ITEM {e*}
 For each $e \in F$, \T e is a binary relation in \Voc, \Function {{\T
e}^{\MC a}} {\Vec \times F} {\Vec \times F}, and for all $(y, b) \in
Vec \times F$,
 \[
	{\T e}^{\MC a} (y, b) = (e \Times_{\Vec} y, e \Times_F b).
 \]
 \end{ITEMS}
\end{DEFINITION}

\Remark
 The cardinality of \Voc is $\aleph_1$ just for the convenience of the
reader. A finite vocabulary is possible by parameterizing the
relations as in \cite [Claim 1.4] {Sh189}.

For each $s \in \Voc \Minus \Set {\Pr [-] n {}} {n < \omega}$, the
interpretation $s^{\MC a}$ is the same for all $\Col a \in
\ColSet$. Hence we omit the superscript \MC a.

For $\mu < \omega_1$, the restriction of \MC a to the set
 \[
	\Par {\Set {y \in \Vec} {\Supp y \Subset \mu} \times F}
	\Union
	\Par {(S \Inter \mu+1) \times \ODFun \omega F}
 \]
 is denoted by \Res {\MC a} {\mu +1}.

\begin{LEMMA} {Equivalent}
 Suppose $\Col a, \Col b \in \ColSet$ and $\mu \leq \omega_1$.
 \begin{ITEMS}

\ITEM {unif->isom}
 If \Function f \mu F uniformizes \Res {(\Col b - \Col a)} {\mu+1},
then $\Res {\MC a} {\mu +1} \Isomorphic \Res {\MC b} {\mu +1}$.

\ITEM {isom->unif}
 If $\Res {\MC a} {\mu +1} \Isomorphic \Res {\MC b} {\mu+ 1}$, then
there is \Function f \mu F which uniformizes \Res {(\Col b - \Col a)}
{\mu +1}.

\ITEM {equivalent}
 $\MC a \LEquiv \MC b$.

 \end{ITEMS}



\begin{Proof}%
\ProofOfItem {unif->isom}
 Suppose \Function f \mu F uniformizes $\Res {(\Col b - \Col a)} {\mu
+1}$. We define \Isomorphism \iota {\Res {\MC a} {\mu +1}} {\Res {\MC
b} {\mu +1}} by the following equations.
 \begin{itemize}

\item
 For all $\xi < \mu$,
 \[
	\iota(\Gen \xi, 0) = \Par {\Gen \xi, f(\xi)},
 \]
 and for all $(y, c) \in \Vec \times F$, we set
 \ARRAY[lll]{
	\iota(y, c)
		&=& \P c \Par [\big] {
			\Sum [\xi < \mu]
			{\T {d_\xi} \Par [\big] {\iota(\Gen \xi, 0)}}
		} \\
		&=& \Par [\Big] {
		\Sum [\xi < \mu] {d_\xi \Gen \xi}, \  
	        \Par [\big] {\Sum [\xi < \mu] {d_\xi \Times f(\xi)}}
		+ c} \\
		&=& \Par [\big] {y, f(y) + c},
 }
 where $y$ is of the form \Sum [\xi < \mu] {d_\xi \Gen \xi}, $d_\xi
\in F$, and $f(y) = \Sum [\xi < \mu] {d_\xi \Times f(\xi)}$ as in
\Section {Preliminaries}.

\item
 For all $\delta \in S \Inter \mu +1$,
 \[
	\iota(\delta, \hat 0) =
	\Par [\big] {\delta, {\hat 0}_\delta^f},
 \]
 where $\hat 0$ denotes the $0$-function of \ODFun \omega F, and
${\hat 0}_\delta^f$ is a function from $\omega$ into $F$ defined for
all $n < \omega$ by
 \ARRAY[lll]{
	{\hat 0}_\delta^f (n)
	&=& \Par [\big] {
		\Sum [\xi < \delta] {e_\xi^{\delta,n} \Times f(\xi)}
	    }
	- (\Col [_{\delta,n}] b - \Col [_{\delta,n}] a) \\
	&=& f(\LadSys [_{\delta,n}])
	- (\Col [_{\delta,n}] b - \Col [_{\delta,n}] a),
 }
 where $\LadSys [_{\delta, n}]$ is of the form \Sum [\xi < \delta]
{e^{\delta, n}_\xi \Times \Gen \xi}, and for all $\xi < \delta$,
$e^{\delta, n}_\xi \in F$.

Furthermore, we define for all $(\delta, u) \in (S \Inter \mu +1)
\times \ODFun \omega F$, that
 \ARRAY[lll]{
	\iota(\delta, u)
	&=& \P u \Par {\iota(\delta, \hat 0)} \\
	&=& \Par {\delta, \  {\hat 0}_\delta^f + u}.
 }

 \end{itemize}
 Since $f$ uniformizes $\Res {(\Col b - \Col a)} {\mu +1}$, the
function ${\hat 0}_\delta^f$ is in \ODFun \omega F for all $\delta \in
S \Inter \mu +1$. Clearly $\iota$ is bijective, and directly by the
definition it preserves all the interpretations of the symbols in
$\Voc \Minus \Set {\Pr [-] n {}} {n < \omega}$. Hence, to prove that
$\iota$ is an isomorphism, it suffices to show that for all $n <
\omega$ and $(\delta, u) \in (S \Inter \mu +1) \times \ODFun \omega
F$,
 \ARRAY[lll]{
 \iota(\Pr n a (\delta, u))
 &=& \iota \Par [\big] {
	\LadSys [_{\delta,n}], \Col [_{\delta,n}] a + u(n)} \\

 &=& \Par [\big] {\LadSys [_{\delta,n}], \ 
	f(\LadSys [_{\delta,n}]) + \Col [_{\delta,n}] a + u(n)} \\

 &=& \Par [\Big] {\LadSys [_{\delta,n}], \ 
	\Col [_{\delta,n}] b
	+ \Par [\big] {f(\LadSys [_{\delta,n}])
		-(\Col [_{\delta,n}] b-\Col [_{\delta,n}] a)}
	+ u(n)} \\

 &=& \Par [\big] {\LadSys [_{\delta,n}], \ 
	\Col [_{\delta,n}] b + {\hat 0}_\delta^f (n) + u(n)} \\

 &=& \Par [\big] {\LadSys [_{\delta,n}], \ 
	\Col [_{\delta,n}] b + ({\hat 0}_\delta^f + u)(n)} \\

 &=& \Pr n b (\delta, \  {\hat 0}_\delta^f + u) \\

 &=& \Pr n b  \Par [\big] {\iota(\delta, u)}.
 }

\ProofOfItem {isom->unif}
 Suppose then \Isomorphism \iota {\Res {\MC a} {\mu +1}} {\Res {\MC b}
{\mu +1}}. We let \Function f \mu F be the unique function satisfying
for all $\xi < \mu$ and $c \in F$, $f(\xi) = c$ iff $\iota(\Gen \xi,
0) = (\Gen \xi, c)$.

Assuming that $\LadSys [_{\delta, n}]$ is of the form \Sum [\xi <
\delta] {e^{\delta, n}_\xi \Times \Gen \xi}, for all $\delta \in S$
and $n < \omega$, the following equation holds in both models,
 \[
	(\LadSys [_{\delta,n}], 0)
	= (\Sum [\xi < \delta] {e^{\delta,n}_\xi \Times \Gen \xi},\; 0)
	= \Sum [\xi < \delta] {e^{\delta,n}_\xi \Times (\Gen \xi,0)}.
 \]
 Hence the isomorphism $\iota$ satisfies
\ARRAY[lll]{
	\iota(\LadSys [_{\delta,n}], 0)
	&=& \Sum [\xi < \delta] {e^{\delta,n}_\xi \Times \iota(\Gen\xi,0)}\\
	&=& \Sum [\xi < \delta] {e^{\delta,n}_\xi \Times (\Gen\xi,f(\xi))}\\
	&=& \Par {\LadSys [_{\delta,n}], f(\LadSys [_{\delta,n}])}.
}
In addition to this, $\iota$ satisfies $\iota(\LadSys [_{\delta, n}],
\Col [_{\delta, n}] a) = \Par {\LadSys [_{\delta, n}], f(\LadSys
[_{\delta,n}]) + \Col [_{\delta, n}] a}$. So the following equation
holds for all $\delta \in S \Inter \mu +1$ and $n < \omega$,
\ARRAY[lll]{
	\Par [\big] {\LadSys [_{\delta,n}],
		f(\LadSys [_{\delta,n}]) + \Col [_{\delta,n}] a}
	&=& \iota(\LadSys [_{\delta,n}], \Col [_{\delta,n}] a) \\
	&=& \iota \Par [\big] {\Pr n a(\delta, \hat 0)} \\
	&=& \Pr n b \Par [\big] {\iota(\delta, \hat 0)} \\
	&=& \Pr n b (\delta, {\hat 0}_\delta^\iota) \\
	&=& (\LadSys [_{\delta,n}],
		\Col [_{\delta,n}] b + {\hat 0}_\delta^\iota(n)),
}
where ${\hat 0}_\delta^\iota$ is the function in \ODFun \omega F
satisfying $\iota(\delta, \hat 0) = (\delta, {\hat 0}_\delta^\iota)$.
It follows that for all $\delta \in S \Inter \mu +1$ and $n < \omega$,
 \[
	\Col [_{\delta,n}] b - \Col [_{\delta,n}] a
	= f(\LadSys [_{\delta,n}]) - {\hat 0}_\delta^\iota(n).
 \]
 Since ${\hat 0}_\delta^\iota \in \ODFun \omega F$, $(\Col b - \Col a)
(\delta) \IsAlmost f(\LadSys (\delta))$ for all $\delta \in S \Inter
\mu +1$, i.e., $f$ uniformizes \Res {(\Col b - \Col a)} {\mu +1}.

\ProofOfItem {equivalent}
 To prove the claim we show that for all $\mu_0 < \mu_1 < \omega_1$
and \Isomorphism {\iota_0} {\Res {\MC a} {\mu_0 +1}} {\Res {\MC b}
{\mu_0 +1}}, there is \Isomorphism {\iota_1} {\Res {\MC a} {\mu_1 +1}}
{\Res {\MC b} {\mu_1 +1}} which is an extension of $\iota_0$. This
suffices by \cite [Theorem 4.3.1 on page 353] {Dickmann}.

By \Item {isom->unif} the existence of $\iota_0$ implies that there is
\Function {f_0} {\mu_0} F uniformizing \Res {(\Col b - \Col a)} {\mu_0
+1}. By \Lemma {Colouring} there is an extension \Function {f_1}
{\mu_1} F of $f_0$ which uniformizes \Res {(\Col b - \Col a)} {\mu_1
+1}. Hence by \Item {unif->isom}, there is \Isomorphism {\iota_1} {\Res
{\MC a} {\mu_1 +1}} {\Res {\MC b} {\mu_1 +1}}.

It can be easily seen from the proof of \Item {isom->unif} that if
$\mu \leq \omega_1$, \Isomorphism {\iota'} {\Res {\MC a} {\mu +1}}
{\Res {\MC b} {\mu +1}}, and \Function f \mu F is the function given
in the proof of \Item {isom->unif}, then the isomorphism $\iota$ given
in the proof of \Item {unif->isom} is the same as $\iota'$. Hence $f_0
\Subset f_1$ implies $\iota_0 \Subset \iota_1$.
\end{Proof}%
%

\end{LEMMA}

\begin{LEMMA}{Isomorphic}
 \begin{ITEMS}

\ITEM {isomorphic}
 For all $a, b \in \ColSet$, $\MC a \Isomorphic \MC b \Iff \Col a
\Equivalent \Col b$.

\ITEM {characterization}
 Suppose \N is a model of vocabulary \Voc, $\Card \N = \aleph_1$, and
$\N \LEquiv \MC a$ for some $\Col a \in \ColSet$. Then there is $\Col
b \in \ColSet$ such that $\N \Isomorphic \MC b$.

\ITEM {card}
  For each $\Col a \in \ColSet$, $\No {\MC a} = \Card {\Quotient
\ColSet \UnifSet}$.

 \end{ITEMS}



\begin{Proof}%
\ProofOfItem {isomorphic}
 The claim holds by \LemmaItem {Equivalent} {unif->isom} and
\LemmaItem {Equivalent} {isom->unif} of \Lemma {Equivalent}.

\ProofOfItem {characterization}
 We let $\phi_\delta$, for all $\delta \in S$, be the following
$\Lan(\Voc)$-sentence,
 \[
 \exists \Seq {r_{\delta,n}} {n < \omega}
 \forall s \in R_\delta
  \Par [\Big] {
   \BigOR [I \in D] {
    \Par [\big] {
     \BigAND [n \in I] {
	\Pr [-] n {}(s) = r_{\delta,n}
     }
     \AND
     \BigAND [n \in \omega \Minus I] {
	\Pr [-] n {}(s) \not= r_{\delta,n}
     }
    }
   }
  }.
 \]
 For all $\delta \in S$, $\phi_\delta$ holds in \N since the
interpretation $r_{\delta, n} = (\LadSys [_{\delta, n}], \Col
[_{\delta, n}] a)$, for all $\delta \in S$ and $n < \omega$, satisfies
the formula in \MC a. We let \Seq {r_{\delta, n}} {n < \omega},
$\delta \in S$, be a sequence of elements in \N satisfying
$\phi_\delta$, and $s_\delta$ be the unique element in ${R_\delta}^\N$
which satisfies ${\Pr [-] n {}}^\N(s_\delta) = r_{\delta, n}$ for all
$n < \omega$.

We define \Function \iota {(S \times \ODFun \omega F) \Union (\Vec
\times F)} \N by the following stipulations.
 \begin{itemize}

\item
 For all $\delta \in S$,
 \[
	\iota(\delta, \hat 0) = s_\delta,
 \]
 \Note {where $\hat 0$ denotes the $0$-function of \ODFun \omega F},
and for all $(\delta, u) \in S \times \ODFun \omega F$,
 \[
	\iota(\delta ,u) = {\P u}^\N \Par {\iota(\delta, \hat 0)}.
 \]

\item
 For all $\xi < \omega_1$, $\iota(\Gen \xi, 0)$ is an arbitrary
element in ${R_{\Gen \xi}}^\N$, and for all $y \in \Vec$,
 \[
	\iota(y ,0)
	= \sideset{}{^\N} \sum_{\xi < \omega_1}
	\Par [\Big] {
	    {(\T {d_\xi})}^\N \Par [\big] {\iota(\Gen \xi, 0)}
	},
 \]
 where $y$ is of the form \Sum [\xi < \mu] {d_\xi \Gen \xi}. For all
$(y, c) \in \Vec \times F$, set $\iota(y, c) = {\P c}^\N \Par
{\iota(y, 0)}$.
 \end{itemize}

Using $\iota$ we define $\Col b$ to be the $F$-colouring on $S$ which
satisfies for all $\delta \in S$ and $n < \omega$,
\[
	\iota(\LadSys [_{\delta,n}], \Col [_{\delta,n}] b)
	= r_{\delta,n}. 
\]
Such a colouring exists since $\iota$ is surjective.

To show that $\iota$ is an isomorphism between \MC b and \N we first
note that $\iota$ is a bijection, and that the preservations of the
interpretations of the symbols in $\Voc \Minus \Set {\Pr [-] n {}} {n
< \omega}$ are obvious. So it suffices to check that $\iota \Par {\Pr
n b (\delta, u)} = {\Pr [-] n {}}^\N \Par {\iota(\delta, u)}$ for all
$n < \omega$ and $(\delta, u) \in S \times \ODFun \omega F$.

For all $u \in \ODFun \omega F$, $n < \omega$, and $s \in
{R_\delta}^\N$,
 \[
	{\P {u(n)}}^\N \Par {{\Pr [-] n {}}^\N(s)}
	=
	{\Pr [-] n {}}^\N \Par {{\P u}^\N(s)},
 \]
 since in \MC a, for all $(\delta, v) \in S \times \ODFun \omega F$,
 \ARRAY[lll]{
	    \P {u(n)} \Par {\Pr n a(\delta,v)}
	&=& \P {u(n)} \Par {
		\LadSys [_{\delta,n}], \  \Col [_{\delta,n}] a + v(n)
	    } \\
	&=& \Par {
		\LadSys [_{\delta,n}], \ 
		\Col [_{\delta,n}] a + v(n) + u(n)
	    } \\
	&=& \Par {
		\LadSys [_{\delta,n}], \ 
		\Col [_{\delta,n}] a + (v+u)(n)
	    } \\
	&=& \Pr n a(\delta, v+u) \\
	&=& \Pr n a \Par {\P u(\delta, v)}.
 }
 Thus for all $n < \omega$ and $(\delta, u) \in S \times \ODFun \omega
F$ the following equation holds,
 \ARRAY[lll]{
	    \iota \Par [\big] {\Pr n b (\delta, u)}
	&=& \iota \Par [\big] {
		\LadSys [_{\delta,n}], \ 
		\Col [_{\delta,n}] b + u(n)
	    } \\
	&=& \iota \Par [\big] {
		\P {u(n)} (\LadSys [_{\delta,n}], \ 
		\Col [_{\delta,n}] b)
	    } \\
	&=& {\P {u(n)}}^\N \Par [\big] {
		\iota(\LadSys [_{\delta,n}], \ 
		\Col [_{\delta,n}] b)
	    } \\
	&=& {\P {u(n)}}^\N (r_{\delta,n}) \\
	&=& {\P {u(n)}}^\N \Par [\big] {
		{\Pr [-] n {}}^\N (s_\delta)
	    } \\
	&=& {\Pr [-] n {}}^\N \Par [\big] {
		{\P u}^\N (s_\delta)
	    } \\
	&=& {\Pr [-] n {}}^\N \Par [\Big] {
		{\P u}^\N \Par [\big] {
			\iota(\delta, \hat 0)
		}
	    } \\
	&=& {\Pr [-] n {}}^\N \Par [\big] {\iota(\delta, u)},
 }
 where we assumed that $\iota$ preserves the interpretations of
symbols \P {u(n)} and \P u.

\ProofOfItem {card}
By \ItemOfLemma {Equivalent} {equivalent} and \Item {isomorphic} \No
{\MC a} is at least \Card {\Quotient {\ColSet [S, F]} {\UnifSet
[\LadSys, D]}}. On the other hand, \Item {characterization} shows that
$\No {\MC a} \leq \Card [\big] {\Set {\Quotient {\MC c} \Isomorphic}
{\Col c \in \ColSet [S, F]}} = \Card {\Quotient {\ColSet [S, F]}
{\UnifSet [\LadSys, D]}}$.
\end{Proof}%
%

\end{LEMMA} 

\begin{PvComment}
 \Remark Generalizations of \Lemma {Equivalent} and \Lemma
{Isomorphic} are straightforward for arbitrary $\kappa > \omega_1$
when $S$ consists of limit ordinals with fixed cofinality and $S
\Inter \mu$ is nonstationary in $\mu$ for all $\mu < \kappa$.
\end{PvComment}


\begin{RefArea}{Env}{Theorem}{Models}

\begin{SeparateProof}
Let $S$ be bistationary in $\omega_1$ and $F$ of size 2. Then by
\Theorem {Colouring} it is consistent with \ZFC + \GCH that there is a
\Vec-ladder system \LadSys on $S$ such that $\Card {\Quotient \ColSet
\UnifSet} = 2$. Then for any $\Col a \in \ColSet$, $\No {\MC a} = 2$
by \ItemOfLemma {Isomorphic} {card}. Now \Theorem {Models} follows
from the following fact \cite {Sh133}:
\begin{property}
 if there is a model \M for which $\No \M = 2$, then for each $k <
\omega$ there is a model $\M_k$ of the same cardinality as \M with
$\No {\M_k} = k$.
 \end{property}
We sketch the proof of this fact. Fix $1 < l < \omega$ and let
$\lambda = \Card \M$. Define $\M_{l +1}$ to be the disjoint union of
$l$-many copies of \M. Add a binary relation symbol $\sim$ to \Voc,
say $\Voc' = \Voc \Union \Braces \sim$, and set for all $x, y \in
\M_{l +1}$ that $x \sim^{\M_{l +1}} y$ iff $x$ and $y$ are in the same
copy of \M. Then each model of cardinality $\lambda$ which is $\Lan
[\lambda] (\Voc')$-equivalent to $\M_{l +1}$ must have the same
structure as $M_{l +1}$ has, i.e., it is a disjoint union of $l$-many
equivalence classes under $\sim$, and each class alone forms a model
$\N_i$, $i < l$, of cardinality $\lambda$ which is $\Lan [\lambda]
(\Voc)$-equivalent to $\M$. Since there are $l +1$-many ways to
select, up to isomorphism, the models $\N_i \LEquiv [\lambda] \M$ for
$i < l$ \Note {the order in the selections of $\N_i$ is immaterial,
only the number of $\N_i$ which are isomorphic to \M matters}, and
because all such selections are pairwise $\Lan [\lambda]
(\Voc')$-equivalent, $\No {\M_{l +1}}$ must be $l +1$.
\end{SeparateProof}

\end{RefArea}


\end{SECTION}



\providecommand{\bysame}{\leavevmode\hbox to3em{\hrulefill}\thinspace}




\catcode`@=11
\begin{tabbing}
Saharon Shelah:\= \\
	\>Institute of Mathematics\\
	\>The Hebrew University\\
	\>Jerusalem. Israel\\
\\
	\>Rutgers University\\
	\>Hill Ctr-Busch\\
	\>New Brunswick. New Jersey 08903\\
	\>\texttt{shelah@math.huji.ac.il}
\end{tabbing}

\begin{tabbing}
Pauli V\"{a}is\"{a}nen:\= \\
	\>Department of Mathematics\\
	\>P.O. Box 4\\
	\>00014 University of Helsinki\\
	\>Finland\\
	\>\texttt{pauli.vaisanen@helsinki.fi}
\end{tabbing}
\catcode`@=13



\end{document}